\documentclass[11pt]{article}
\pdfoutput=1
\usepackage{amssymb,amsmath,amsfonts,mathrsfs,bm,color,galois,graphicx}
\usepackage{latexsym}
\usepackage{booktabs}
\usepackage{subfigure}
\usepackage{multirow,makecell}
\usepackage{float}
\usepackage[ruled,linesnumbered]{algorithm2e}
\SetKwRepeat{Do}{do}{while}
\usepackage[figuresright]{rotating}
\usepackage[toc,page]{appendix}
\usepackage[round,authoryear]{natbib}
\usepackage[colorlinks, linkcolor=blue, anchorcolor=blue, citecolor=blue]{hyperref}
\usepackage{enumitem}

\usepackage{array}
\renewcommand{\d}{\mathrm{d}}

\allowdisplaybreaks
 \numberwithin{equation}{section}
 \numberwithin{theorem}{section}

\title{A Less Uncertain Sampling-Based Method\\of Batch Bayesian Optimization}
\author{{Kai Jia$^{1}$ \quad Xiaojun Duan$^{1}$\footnote{Corresponding author. Email: xjduan@nudt.edu.cn} \quad Zhengming Wang$^{2}$ \quad Liang Yan$^{1}$}\\[6pt]
 { 1. College of arts and sciences, National University of Defense Technology}\\
 { 2. College of systems engineering, National University of Defense Technology}\\
 { Changsha, Hunan, 410005, China} \\[8pt]}
\date{\today}

\begin{document}
\baselineskip=18pt
\maketitle

\begin{abstract}  
  This paper presents a method called sampling-computation-optimization (SCO) to design batch Bayesian optimization.
  SCO does not construct new high-dimensional acquisition functions but samples from the existing \textit{one-site} acquisition function to obtain several candidate samples.
  To reduce the uncertainty of the sampling, the general discrepancy is computed to compare these samples.
  Finally, the genetic algorithm and switch algorithm are used to optimize the design.
  Several strategies are used to reduce the computational burden in the SCO.
  From the numerical results, the SCO designs were less uncertain than those of other sampling-based methods.
  As for application in batch Bayesian optimization, SCO can find a better solution when compared with other batch methods in the same dimension and batch size.
  In addition, it is also flexible and can be adapted to different one-site methods.
  Finally, a complex experimental case is given to illustrate the application value and scenario of SCO method.
\bigskip

\noindent{\bf  Keywords}: 
  Batch sequential design; Rejection sampling; Sampling-importance-resampling; General discrepancy; Global optimization
\end{abstract}

\section{Introduction}

Bayesian optimization is important in the engineering optimization field.
The general idea of Bayesian optimization is as follows:
\begin{enumerate}
	\item[(i)] Build or update the surrogate model $\hat{f}_{\rm exp}$ based on the database $Y_{\rm data}=f_{\rm exp}(X_{\rm data})$; $\hat{f}_{\rm exp}$ usually provides the prediction or description of the experiment response and is easy to solve;
	
	\item[(ii)] Construct the acquisition function $\phi(\bm x)=\phi(\bm x|X_{\rm data},Y_{\rm data},\hat{f}_{\rm exp})$ to evaluate the value of further experiments on $\Omega$ according to the model and experiment purposes;
	
	\item[(iii)] Choose the most valuable sites $\bm x^*=\arg\max\phi(\bm x)$ by some auxiliary optimization to arrange the sequential experiments;
	
	\item[(iv)] Conduct the sequential experiments, and update the database.
\end{enumerate}

The acquisition function $\phi$ is usually the function of one site.
The classic Bayesian efficient global optimization (EGO) method \citep{jones1998efficient} uses a  proposed expected improvement (EI), based on a Gaussian process model as an acquisition function, which is efficient and widely used.
According to (iii), the sequential experiments only include one site in each experiment period.
However, this is very time-consuming, especially when there are many parallel experimental resources.
In order to make better use of these resources, the experimenter needs to identify multiple sequential experimental sites for each period.

One way to implement this is to choose several sub-optimal sites under $\phi$, but these sub-optimal sites may converge to a cluster, which make the experiment inefficient.
The other way is to construct a multi-point acquisition function $\phi(\bm x_1,\dots,\bm x_n)$, such as $q$-EI \citep{ginsbourger2008a,ginsbourger2010kriging}.
However, a multi-point acquisition function is usually highly dimensional, which make it difficult to compute and optimize.

To solve the above problems, sampling technique is used to implement batch Bayesian optimization.
\cite{cai2017a} and \cite{wang2004mode} use the mode-pursuing sampling method \citep{fu2002a} to obtain the sequential design.
\cite{ning2020batch} use the sampling-importance-resampling method to accelerate EGO.
However, sampling methods are arbitrary and uncertain, which is not appropriate in an expensive experiment.

This article presents a method called sampling-computation-optimization (SCO) to design the sequential experiments of batch Bayesian optimization.
SCO does not construct new high-dimensional acquisition functions; instead it samples from the existing one-site acquisition function to obtain several candidate samples.
Some calculations and optimizations then proceed to obtain the sequential design.
The rest of this paper is organized as follows.
Section \ref{sec pre} introduces two important sampling methods and the general discrepancy concept.
The details of the SCO method are described in Section \ref{sec SCO}.
Some numerical analyses and comparisons are presented in Section \ref{sec num}.

\section{Preliminary}\label{sec pre}

\subsection{Rejection sampling and sampling-importance-resampling}

Rejection sampling (RS) and sampling-importance-resampling (SIR) are two widely used sampling methods.
We introduce the two methods briefly, for more information, see \cite{bishop2006pattern}.

Suppose we would like to sample from a target density $f(\cdot)\propto\phi(\cdot)$ on $\Omega$, which is hard to sample directly but easy to compute.
The main steps of RS are as follows:
\begin{enumerate}[itemindent=2.3em]
	\item[Step 0]: Choose a proposal distribution with density $g(\cdot)$, which is easy to sample.
	
	\item[Step 1]: Find a constant $M$ that is large enough to ensure $\phi(\bm x)\le Mg(\bm x),\ \forall \bm x\in\Omega$.
	
	\item[Step 2]: Generate $\bm u_i\sim g(\cdot)$ and $v_i\sim {\rm U}(0,1)$.
	
	\item[Step 3]: If $\phi(\bm u_i)<Mg(\bm u_i)v_i$, then accept $\bm u_i$ as a sample $\bm x_j$.
	Otherwise, reject $\bm u_i$.
\end{enumerate}
As indicated by \cite{bishop2006pattern}, $\bm x_j\sim f(\cdot)$, and the probability to accept $\bm u_i$ is inversely proportional to $M$.
Generally, one needs to choose a smallest $M$ to ensure $\phi(\bm x)\le Mg(\bm x),\ \forall \bm x\in\Omega$.

Sometimes, a suitable $M$ is difficult to find, or the reject rate is too high.
We then need to turn to SIR.
The main steps of SIR are as follows:
\begin{enumerate}[itemindent=2.3em]
	\item[Step 0]: Choose a proposal distribution with density $g(\cdot)$, which is easy to sample.
	
	\item[Step 1]: Generate $U=\{\bm u_1,\dots,\bm u_N\},\ \bm u_i\sim g(\cdot)$.
	
	\item[Step 2]: Calculate the weights $w_i=\frac{\phi(\bm u_i)/g(\bm u_i)}{\sum_{i=1}^N\phi(\bm u_i)/g(\bm u_i)}$.
	
	\item[Step 3]: Resample $\bm x_j$ from $U$ with the probability $(w_1,\dots,w_N)$.
\end{enumerate}
\cite{bishop2006pattern} indicated that $\bm x_j\sim f(\cdot)$ when $N\rightarrow\infty$, which means SIR needs a large amount of pre-samples.

Both RS and SIR first need to sample $\bm u_i$ from the proposal distribution.
For convenience, we call this process \emph{pre-sampling} and refer to U as a \emph{pre-sample set}.
The two sampling methods are respectively suitable for two situations of the algorithm in this paper, and we will introduce their application in Section 3.

\subsection{General Discrepancy}

There are many interpretations of the discrepancy \citep{li2020is}.
In this paper, we refer to the interpretation in Section 2.4 of \cite{fang2018theory}.
We will begin with some basic notation.

Let $\Omega$ represent the experimental domain, and we only discuss the case when $\Omega=[0,1]^d$.
The design is taken as a set $X=\{\bm x_1,\bm x_2,\dots,\bm x_n\}$ on $\Omega$.
$K(\cdot,\cdot):\Omega\times\Omega\rightarrow\mathcal{R}$ is a symmetric and positive semi-definite kernel function.
$F$ is a distribution function in 
$$\mathcal{K}=\left\{F\left|\int_{\Omega\times\Omega}K(\bm u,\bm v)\d F(\bm u)\d F(\bm v)<\infty\right.\right\}.$$
Define the inner product of two arbitrary functions $F,G\in\mathcal{K}$ as
$$\left<F,G\right>_K=\int_{\Omega\times\Omega}K(\bm u,\bm v)\d F(\bm u)\d G(\bm v).$$
The general discrepancy of the design $X$ with respect to the target distribution $F$ using the kernel $K$ is then defined as
\begin{align}\label{eq D2}
D^2(X,F,K)\triangleq&||F_X-F||_K^2=\left<F_X-F,F_X-F\right>_K \nonumber \\
=&\int_{\Omega\times\Omega}K(\bm{u},\bm{v})\d (F-F_X)(\bm{u})\d (F-F_X)(\bm{v}) \nonumber \\ 
=&\int_{\Omega\times\Omega}K(\bm{u},\bm{v})\d F(\bm{u})\d F(\bm{v}) \nonumber \\
&-\frac{2}{n}\sum_{i=1}^n\int_{\Omega}K(\bm{u},\bm{x}_i)\d F(\bm{u})
+\frac{1}{n^2}\sum_{i,j=1}^nK(\bm{x}_j,\bm{x}_i).
\end{align}
This can be taken as the distance between the target distribution $F$ and the empirical distribution of $X$.
A low discrepancy means $F_X$ and $F$ are close in some sense, so we prefer a design with low discrepancy to better represent the target distribution.

The kernel function $K$ corresponds to an inner product of the Hilbert space.
By taking different functions, we can define different kinds of discrepancy, such as the widely used centered discrepancy (CD), wrapped discrepancy (WD), or mixed discrepancy (MD).
In this article, we use the kernel function of the WD \citep{fang2018theory}:
\begin{equation}\label{eq KW}
K(\bm{u},\bm{v})=\prod_{i=1}^d\left[\frac{3}{2}-|u_i-v_i|+(u_i-v_i)^2\right].
\end{equation}

Because $F$ is arbitrary, the high-dimension integral in equation \eqref{eq D2} is difficult to compute; we use the Monte Carlo method to estimate $D^2(X,F,K)$.
Rewrite equation \eqref{eq D2} as
\begin{align}\label{eq D2 r}
D^2(X,F,K)=&\int_{\Omega\times\Omega}K(\bm{u},\bm{v})\cdot f(\bm{u})\cdot f(\bm{v})\d\bm{u}\d\bm{v} \nonumber \\
&-\frac{2}{n}\sum_{i=1}^n\int_{\Omega}K(\bm{u},\bm{x}_i)\cdot f(\bm{u})\d\bm{u}
+\frac{1}{n^2}\sum_{i,j=1}^nK(\bm{x}_j,\bm{x}_i).
\end{align}
We then get an estimation of $D^2(X,F,K)$ as
\begin{align}\label{eq D2 es}
\hat{D}^2(X,F,K)=&\frac{1}{N^2}\sum_{i,j=1}^NK(\bm{u}_j,\bm{u}_i)\cdot f(\bm{u}_j)\cdot f(\bm{u}_i)\nonumber \\
&-\frac{2}{nN}\sum_{i=1}^n\sum_{j=1}^NK(\bm{u}_j,\bm{x}_i)\cdot f(\bm{u}_j)
+\frac{1}{n^2}\sum_{i,j=1}^nK(\bm{x}_j,\bm{x}_i),
\end{align}
where $\bm{u}_i\sim {\rm U}(\Omega), i=1,2,\dots,N$.

\section{Method for batch Bayesian optimization}\label{sec SCO}

In this paper, we propose a method to implement batch Bayesian optimization.
The new method is based on a one-site acquisition function, and its main idea is to exchange (iii) in Section 1 with a revised (iii') as follows:
\begin{enumerate}
	\item[(iii')] Take the acquisition function $\phi(\bm x)$ as the density of a distribution $F$ and construct a sequential design $X^*$ to fit $F$ best;
\end{enumerate}

For convenience, note $D^2(X,\phi,K)\triangleq D^2(X,F,K)$, where $F$ is induced by density $f\propto\phi$.
In addition, we use $D^2(X,\phi,K)$ to measure the fitness described above in (iii').
If $F\sim {\rm U}(\Omega)$, this kind of sequential design is also called a uniform design.
The methods for constructing uniform design are various and efficient, such as number theorem methods \citep{niederreiter1992random,fang1994number,fang1994some}, algorithmic optimization methods \citep{fang2000uniform,zhou2012constructing,zhou2013an,chen2014discrete}, and hybrid methods \citep{zhou2013mixture}.

However, when F is arbitrary, it becomes complicated.
The first difficulty is that number theorem methods are not suitable when $F$ is not uniform.
The second is $D^2(X,F,K)$ has no analytic form when $F$ is arbitrary.
Calculating $D^2(X,F,K)$ is already difficult, let alone optimizing it.
As a result, the construction has to make a trade-off between optimality and feasibility, aiming to construct a design with relatively low discrepancy in an acceptable time.

In this section, we first used the sampling method to generate several candidate samples.
To avoid the arbitrariness of sampling method, general discrepancy is calculated to compare candidate designs.
Ultimately, the sequential design is determined through optimization.
We call this method sampling-computation-optimization (SCO).
Next, the details of SCO are described step-by-step.

\subsection{Sampling}

At the beginning, we use RS to obtain samples.
The proposal distribution $g(\cdot)$ is chosen as ${\rm U}(\Omega)$ for convenience.
In this way, the pre-sample set $U$ can be used in calculation, and we will explain this next.
As a requirement of RS, $M$ should be large enough to ensure $\phi(\bm x)\le Mg(\bm x)=M$, we need to optimize $\phi(\bm x)$ to determine the minimum of $M$.
The advantage in doing this is so that we can find the optimal $\bm x^*$ and take it as the first site of $X$; we then only sample $n-1$ points to form the candidate design.
In this way, the batch Bayesian optimization can be compatible with the one-site Bayesian optimization.
As mentioned previously, the probability of accepting samples in RS is inversely proportional to $M$.

\begin{figure}[htb]
	\centering
	\subfigure[flat at the beginning]{
		\begin{minipage}[t]{0.45\textwidth}
			\centering
			\includegraphics[width=\textwidth]{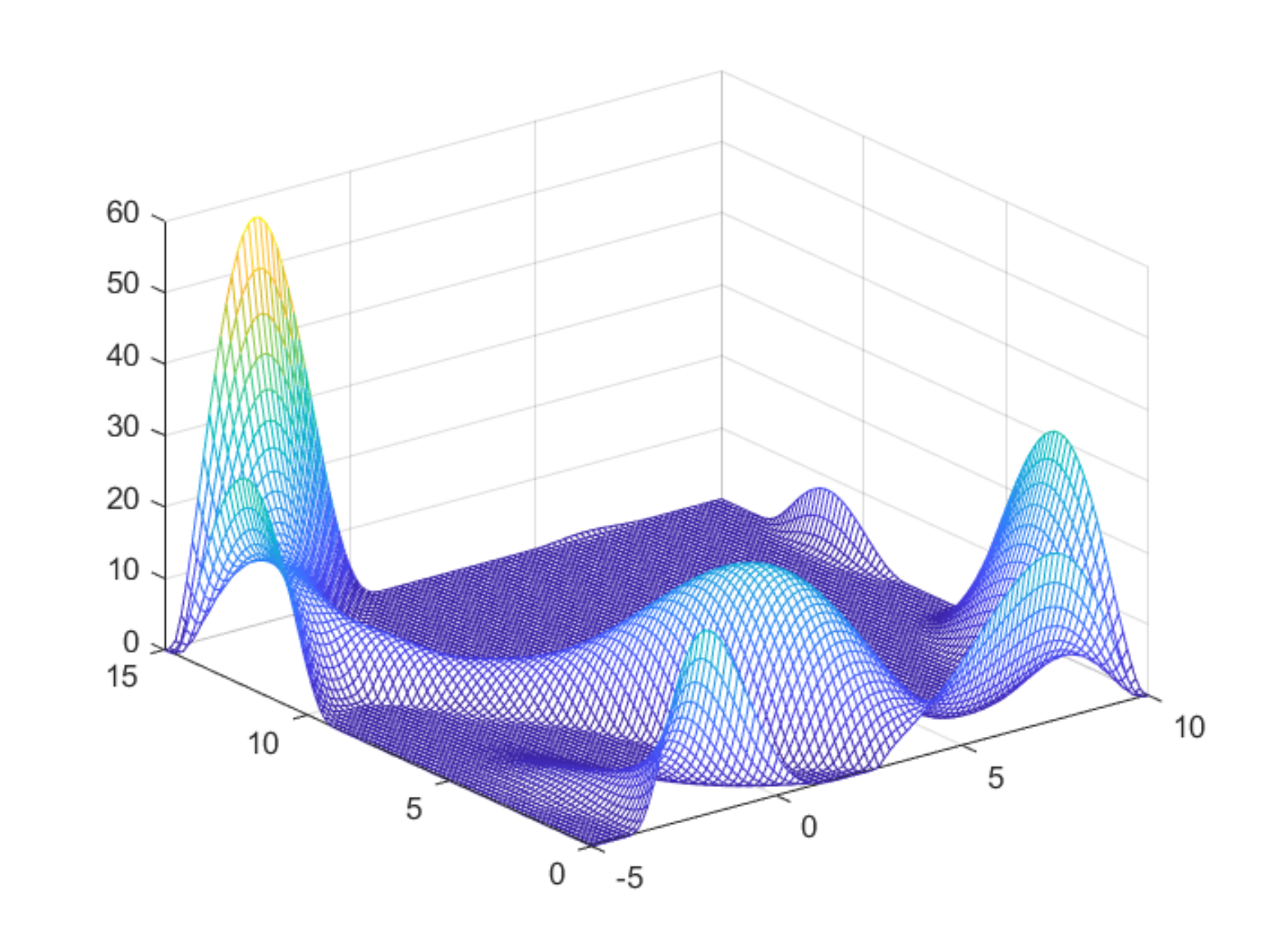}
		\end{minipage}%
	}
	\centering
	\subfigure[steep at the end]{
		\begin{minipage}[t]{0.45\textwidth}
			\centering
			\includegraphics[width=\textwidth]{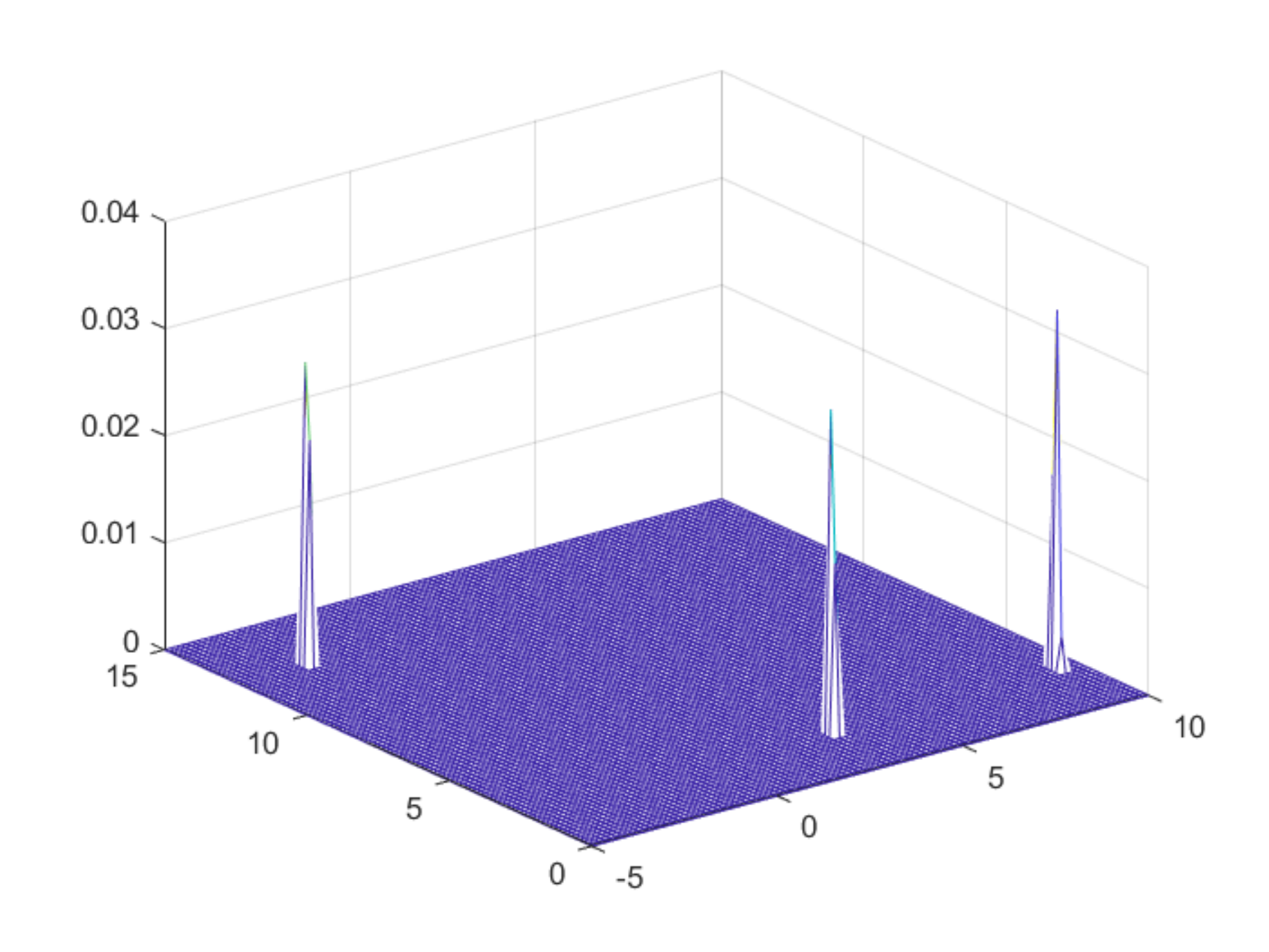}
		\end{minipage}%
	}
	\caption{Two types of EI in EGO}
	\label{fig: flat}
\end{figure}

Considering the rejection rate, there are two drawbacks as follows.
First, when $\phi$ is relatively flat, the rejection rate is too low, and the size of $U$ is too small.
Because $U$ will be used in the calculation, this will affect the accuracy of the next stage.
Second, when $\phi$ is steep, the rejection rate is too high, and the sampling efficiency will be too low, resulting in a large amount of calculation.
These two situations are common at the beginning and end of the EGO algorithm respectively, for example, see Figure \ref{fig: flat}.
We then improve the sampling process to suit the above situations.

For the first drawback, we generate $N_{\rm min}$ points $\bm{u}_i$ and $v_i$ at once, where $N_{\rm min}$ is set to ensure the accuracy of the calculation.
Instead of setting a fixed $M$, the rejection factor is defined by $\lambda_i=v_i\cdot \phi_{\rm max}/\phi(\bm u_i)$.
If the $(n-1)$-th smallest $\lambda_{(n-1)}\le1$, there are at least $n-1$ points accepted.
We then sub-sample from the accepted samples by choosing the $n-1$ points with the smallest $\lambda_i$.
This is equivalent to setting $M=\phi_{\rm max}/\lambda_{(n-1)}$.
If $\lambda_{(n-1)}>1$, we set $M=\phi_{\rm max}$ and proceed with the rejection sampling.

To justify this, assume the target distribution is $f(\bm x)$, the proposal distribution is $g(\bm x)$, and $Mg(\bm x)\ge f(\bm x),\forall \bm x$, $\lambda=\frac{Mg(\bm x)}{f(\bm x)}v, v\sim {\rm U}(0,1)$. Thus
\begin{align}\label{eq dis of lambda}
\pi(\lambda|\textbf{accept},\bm x)=&\frac{\pi(\textbf{accept}|\lambda,\bm x)\cdot\pi(\lambda|\bm x)}{\pi(\textbf{accept}|\bm x)} \nonumber \\
=&\frac{I_{0\le\lambda\le1}\cdot f(\bm x)/(Mg(\bm x))}{f(\bm x)/(Mg(\bm x))} \nonumber \\
=&1\cdot I_{0\le\lambda\le1}\sim {\rm U}(0,1).
\end{align}
This means the distribution of $\lambda$, under the condition of acceptance, is ${\rm U}(0,1)$, which is independent of $\bm x$.
The sub-samples according $\lambda$ then have the same distribution of the original samples.

For the second drawback, if we continue the above RS until $|U|\ge N_{\rm max}$ and the accepted samples are still not enough, we then turn to SIR, and $U$ in RS is directly regard as the pre-sample in SIR.
$N_{\rm max}$ is set to accommodate the computational burden, and it is large enough to meet the requirements of SIR.

\subsection{Calculation}

The sampling and calculation can be carried out simultaneously to reduce the computational burden.
Upon reviewing equation \eqref{eq D2 es}, if we take the proposal distribution $g(\cdot)$ as ${\rm U}(\Omega)$, then the pre-sample set $U$ in RS and SIR can be used to estimate $D^2(X,\phi,K)$.
RS and SIR will be suitable for the calculation of $\hat{D}^2(X,\phi,K)$, and only minor modifications are required.

To take an acquisition function as a density, we need to normalize it by $f(\bm x)=\phi(\bm x)/\int_\Omega\phi(\bm u)\d\bm u$, where $\int_\Omega\phi(\bm u)\d\bm u$ can be estimated by
\begin{equation}
\int_\Omega\phi(\bm u)\d\bm u\ \hat{=}\ \frac{1}{N}\sum_{i=1}^N\phi(\bm u_i)\triangleq\frac{S_\phi}{N}.
\end{equation}
Referring then to equation \eqref{eq D2 es}, $\hat{D}^2(X,\phi,K)$ can be written as
\begin{align}\label{eq D2 es c}
\hat{D}^2(X,\phi,K)=&\frac{1}{S_\phi^2}\sum_{i,j=1}^NK(\bm{u}_j,\bm{u}_i)\cdot \phi(\bm{u}_j)\cdot \phi(\bm{u}_i)\nonumber \\
&-\frac{2}{nS_\phi}\sum_{i=1}^n\sum_{j=1}^NK(\bm{u}_j,\bm{x}_i)\cdot \phi(\bm{u}_j)
+\frac{1}{n^2}\sum_{i,j=1}^nK(\bm{x}_j,\bm{x}_i).
\end{align}

As mentioned previously, sampling once without comparison is arbitrary and uncertain.
We would like to generate several sets of samples as candidates, noted as $X^{(1)},X^{(2)},\dots,X^{(m)}$.
Because the calculation above is very computationally intensive, we propose three measures to reduce this:

\begin{enumerate}
	\item Reduce the calculation in pre-sampling.
	Every time we sample, keep $U$ unchanged and generate $X^{(k)}$.
	For RS, we only regenerate $v_i\sim {\rm U}(0,1)$ to screen $U$ again;
	For SIR, we resample directly from $U$.
	In this way, we only calculate $\phi(\bm u_i)$ once, no matter how many times we sample.

	\item Eliminate unnecessary calculations in $\hat{D}^2(X,\phi,K)$.
	Rewrite equation \eqref{eq D2 es c} as
	\begin{align}\label{eq D2 es c S}
	\hat{D}^2(X,\phi,K)=&A_1(\phi,K)-\frac{2}{n}\sum_{i=1}^nA_2(\bm x_i,\phi,K)+A_3(X,K),
	\end{align}
	where $A_1(\phi,K)=\sum_{i,j=1}^NK(\bm{u}_j,\bm{u}_i)\cdot \phi(\bm{u}_j)\cdot \phi(\bm{u}_i)/S_\phi^{2}$ is irrelative to $X$, $A_3(X,K)=\sum_{i,j=1}^n K(\bm{x}_j,\bm{x}_i)/n^{2}$ is irrelative to $\phi$, and $A_2(\bm x_i,\phi,K)=\sum_{j=1}^NK(\bm{u}_j,\bm{x}_i)\cdot \phi(\bm{u}_j)/S_\phi$ is irrelative to other elements of $X$.
	Because $A_1$ is independent of $X$ and it is most computationally intensive, define
	\begin{align}\label{eq D2 es c minus}
	\hat{D}_-^2(X,\phi,K)=-\frac{2}{n}\sum_{i=1}^nA_2(\bm x_i,\phi,K)+A_3(X,K).
	\end{align}
	We only compute equation \eqref{eq D2 es c minus} to compare the samples.
	If an exact $\hat{D}^2(X,\phi,K)$ is necessary for analysis, we calculate $A_1(\phi,K)$ once at the end, and add it to $\hat{D}_-^2(X,\phi,K)$.
	
	\item Save the important calculation results.
	In addition to the candidate samples $X^{(k)}$ and their corresponding general discrepancies $D^{(k)}$, we save other calculation results for future optimization.
	Firstly, we save the sample set as $S=\cup_{k=1}^mX^{(k)}$.
	Next, because $A_2$ is only related to one site, $A_2(\bm s_i,\phi,K)$ are saved as $A(i)$ for future optimization, where $\bm s_i$ are the corresponding elements of $S$.
	In addition, $\phi(\bm{u}_i)$ calculated in the sampling process are saved as $\Phi(i)$ and are retrieved when calculating $K(\bm{u}_j,\bm{x}_i)\phi(\bm{u}_i)$ in $A_2(\bm s_i,\phi,K)$.
	
\end{enumerate}
The complete processes of sampling and calculation are summarized in Algorithm \ref{SC}.

\begin{algorithm}\label{SC}
	\caption{Sample and Calculation of SCO}
	\KwIn{One-site acquisition function $\phi$, design size $n$, minimum sample size $N_{\rm min}$, maximum sample size $N_{\rm max}$, and number of candidate samples $m$}
	\KwOut{Candidate samples $X^{(k)}$, sample set $S$, calculation results $D^{(k)}$, and $A$}
	let $N=N_{\rm min}$ and generate $U=\{\bm u_i,i=1,\dots,N\}$ uniformly distributed on $\Omega$\;
	calculate $\phi_i=\phi(\bm u_i),i=1,\dots,N$ and save them as $\Phi$\;
	find $\phi_{\rm max}=\max \phi(\bm x)=\phi(\bm x^*)$ by some global optimization algorithm\;
	initialize $S=\{\bm x^*\},\ A=\{A_2(\bm x^*,\phi,K)\}$\;
	\For{$k=1:m$}{
		$X^{(k)}=\{\bm x^*\}$\tcp*{maintain $\bm x^*\in X$}
		\eIf(\tcp*[f]{RS}){$N<N_{\rm max}$}{
			generate $v_i,i=1,\dots,N$ uniformly distributed on $(0,1)$\;
			calculate $\lambda_i=(v_i\cdot \phi_{\rm max})/\phi_i,i=1,\dots,N$\;
			find the $(n-1)$-th smallest $\lambda_i=\lambda_{(n-1)}$\;
			\eIf(\tcp*[f]{Subsampling from RS}){$\lambda_{(n-1)}\le1$}{
				$X^{(k)}=X^{(k)}\cup\{\bm u_i|\lambda_i\le\lambda_{(n-1)}\}$\;
			}(\tcp*[f]{Continue to RS}){
				$X^{(k)}=X^{(k)}\cup\{\bm u_i|\lambda_i\le1\}$\;
				\While{$|X^{(k)}|<n$}{
					$N=N+1$\;
					generate $\bm u_N\sim {\rm U}(\Omega)$, $U=U\cup\{\bm u_N\},\ \Phi=\Phi\cup\{\phi_N=\phi(\bm u_N)\}$\;
					generate $v_N\sim {\rm U}(0,1)$\;
					\If{$v_N\cdot \phi_{\rm max}\le \phi_N$}{
						$X^{(k)}=X^{(k)}\cup\{\bm u_N\}$\;
					}
					\If(\tcp*[f]{Prepare for SIR}){$N\ge N_{\rm max}$}{
						calculate $w_i=\phi_i/\sum_{j=1}^{N}\phi_j,i=1,\dots,N$\;
						break and turn to SIR\;
					}
				}
			}
		}(\tcp*[f]{SIR}){
			resample $\bm x^{(k)}_2,\dots,\bm x^{(k)}_n$  from $U$ with the probability $(w_1,\dots,w_N)$\;
			$X^{(k)}=X^{(k)}\cup\{\bm x^{(k)}_2,\dots,\bm x^{(k)}_n\}$\;
		}
		$S=S\cup X^{(k)}$\;
		\For{$i=2:n$}{
			calculate $A_2(\bm x^{(k)}_i,\phi,K)$ by $\bm x^{(k)}_i,U$ and $\Phi$\;
			$A=A\cup\{A_2(\bm x^{(k)}_i,\phi,K)\}$\;
		}
		calculate $D^{(k)}=\hat{D}_-^2(X^{(k)},\phi,K)$ by $X^{(k)}$ and $A_2(\bm x^{(k)}_i,\phi,K)$\;
	}
	calculate $\hat{D}_+=A_1(\phi,K)$ by $U$ and $\Phi$\tcp*{If exact $\hat{D}$ is necessary}
	For all $k$, $D^{(k)}=D^{(k)}+\hat{D}_+$\;
\end{algorithm}

\subsection{Optimization}\label{SCO-O}

After sampling and calculation, $X^{(1)},X^{(2)},\dots,X^{(m)}$ are still uncertain.
To reduce the uncertainty of the final design $X^*$, we would like to optimize them under the criterion of general discrepancy.
There are two difficulties in optimization.
One is the large number of combinations of candidate sets; if we choose $X$ from $U$, all possible combinations are $\binom{N}{n}$, which makes optimization impractical.
The other difficulty is the burden of calculation.
No matter what kind of algorithm we choose, every time we generate a new design $X^{\rm new}$, calculating $\hat{D}^2(X^{\rm new},\phi,K)$, even $\hat{D}_-^2(X^{\rm new},\phi,K)$, is time-consuming.

Therefore, we need to optimize on a relative-small but reasonable set and make full use of the existing calculation results.
Referring to equation \eqref{eq D2 es c S}, if $X^{\rm new}$ is generated from an old design $X^{\rm old}$ with general discrepancy $D^{\rm old}$ and all the sites are restricted in $S$, then $D^{\rm new}$ can be updated by
\begin{align}\label{eq update}
D^{\rm new}=&D^{\rm old}+\frac{2}{n}\sum_{i\in I_{\rm old}}A(i)-\frac{2}{n}\sum_{i\in I_{\rm new}}A(i)-A_3(X^{\rm old},K)+A_3(X^{\rm new},K) \nonumber \\
\triangleq&D^{\rm old}+\Delta(X^{\rm new},X^{\rm old}),
\end{align}
where $D^{\rm old}$ and $A(i)$ are already known.
We only calculate $A_3(X^{\rm old},K)$ and $A_3(X^{\rm new},K)$, and retrieve the changed terms of $A(i)$ in $I_{\rm old},I_{\rm new}$ to update $D$.

In the rest of this subsection, we introduce two algorithms to optimize $X^*$.


\noindent\textbf{The genetic algorithm}

The genetic algorithm (GA) is a classical intelligent optimization algorithm.
We briefly introduce the process of the GA in Algorithm \ref{GA} and explain the operations of some important steps.
For more details of the GA, one can refer to \cite{goldberg1988genetic}.

\begin{algorithm}\label{GA}
	\caption{Genetic algorithm}
	\KwIn{Data: candidate designs $X^{(k)}$, sample set $S$, calculation results $D^{(k)}$, and $A$;\\
	Parameters: $l$ in fitness function, crossover probability $P_{\rm c}$, and mutation probability $P_{\rm m}$}
	\KwOut{Optimal design $X^*$ and its corresponding $D^*$}
	initialize population: $Parents=\{X^{(1)},\dots,X^{(m)}\}$\;
	$k^*=\arg\min_kD^{(k)},X^*=X^{(k^*)},D^*=D^{(k^*)}$\;
	\While{Stopping criterion is not met}{
		$Offspring=\{X^*\}$\tcp*{Retain the best individual}
		calculate the fitness function ${\rm Fit}(X)=D^{-l}(X)$\;
		\While{$|Offspring|<m$}{
			$(Father,Mother)=Selection(Parents;{\rm Fit})$\;
			$Child=Crossover(Father,Mother;P_{\rm c})$\;
			$Child=Mutation(Child,S;P_{\rm m})$\;
			$D^{Child}=D^{Father}+\Delta(D^{Child},D^{Father})$\;
			$Offspring=Offspring\cup\{Child\}$\;
		}
		update $X^*$ and $D^*$\;
		$Parents=Offspring$\;
	}
\end{algorithm}

\begin{itemize}
\item Fitness function---${\rm Fit}(X)=D^{-l}(X)$;
\item Selection---two parents were selected using the roulette wheel method according to the ${\rm Fit}$ function; the best individual in each iteration remained unchanged to the next generation;
\item Crossover---offspring were generated by exchanging the corresponding sites in parents with probability $P_{\rm c}$;
\item Mutation---replace sites in offspring with a random site in $S$ with probability $P_{\rm m}$, except the first site $\bm x^*$.
\end{itemize}

The default parameters of the GA are set as $l=5,\ P_{\rm c}=0.5,\ P_{\rm m}=0.1$.

\noindent\textbf{The switching algorithm}

The switching algorithm (SA) is an efficient but local algorithm.
In the design literature, it has been widely applied to fit the large space of optimization \citep{winker1998optimal,2000Centered,chuang2010uniform}.
The main steps of the SA are shown in Algorithm \ref{SA}, and it requires no parameters.

\begin{algorithm}\label{SA}
	\caption{Switch algorithm}
	\KwIn{Candidate designs $X^{(k)}$, sample set $S$, calculation results $D^{(k)}$, and $A$}
	\KwOut{Optimal design $X^*$ and its corresponding $D^*$}
	$k^*=\arg\min_kD^{(k)},X^*=X^{(k^*)},D^*=D^{(k^*)}$\;
	\Do{$X^*$ has changed}{
		\For(\tcp*[f]{Maintain $\bm x^*\in X^*$}){i=2:n}{
			find $\bm s^*=\arg\min_{\bm s\in S}\Delta(X_{i,\bm s}^{\rm switch},X^*)$,
			where $X_{i,\bm s}^{\rm switch}=X^*\backslash\{\bm x_i\}\cup\{\bm s\}$\;
			\If{$\Delta^*=\Delta(X_{i,\bm s^*}^{\rm switch},X^*)<0$}{
				$X^*=X_{i,\bm s^*}^{\rm switch}$\;
				$D^*=D^*+\Delta^*$\;
			}
		}
	}
\end{algorithm}

From Algorithm \ref{GA} and \ref{SA}, we can see both the GA and SA have the following properties:

\begin{itemize}
	\item $\bm x^*$ is always the first site of $X^*$;
	\item $D^*$ is monotonically decreasing;
	\item every site in iterations is restricted in $S$;
	\item the inputs of Algorithm \ref{GA} and \ref{SA} are the outputs of Algorithm \ref{SC}.
\end{itemize}

\subsection{Summary}

\begin{figure}[htbp]
	\centering
	\includegraphics[scale=0.6]{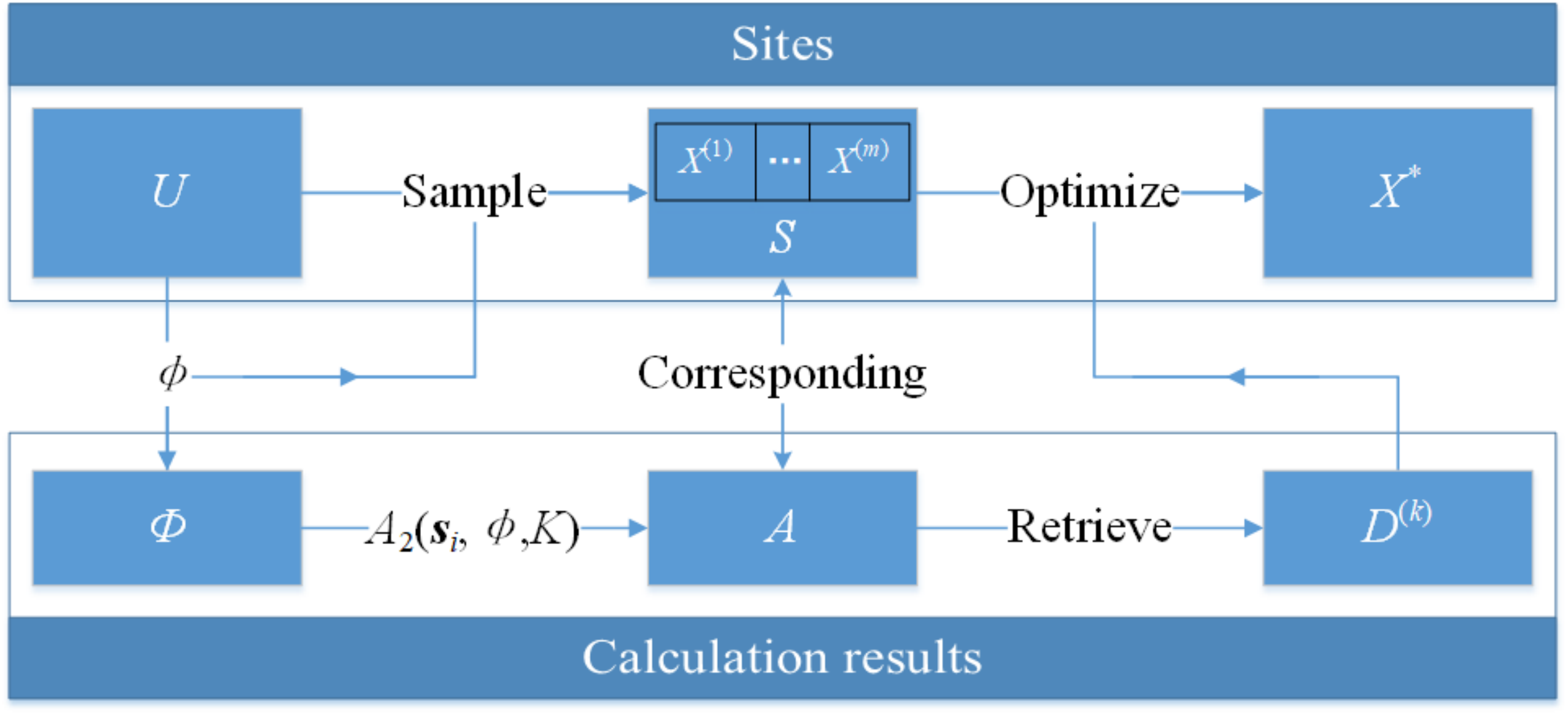}
	\caption{Flowchart of SCO}
	\label{fig: flow SCO}
\end{figure}

\begin{figure}[htb]
	\centering
	\subfigure[Plots of different sets in SCO]{
		\begin{minipage}[t]{.45\textwidth}
			\centering
			\includegraphics[width=\textwidth]{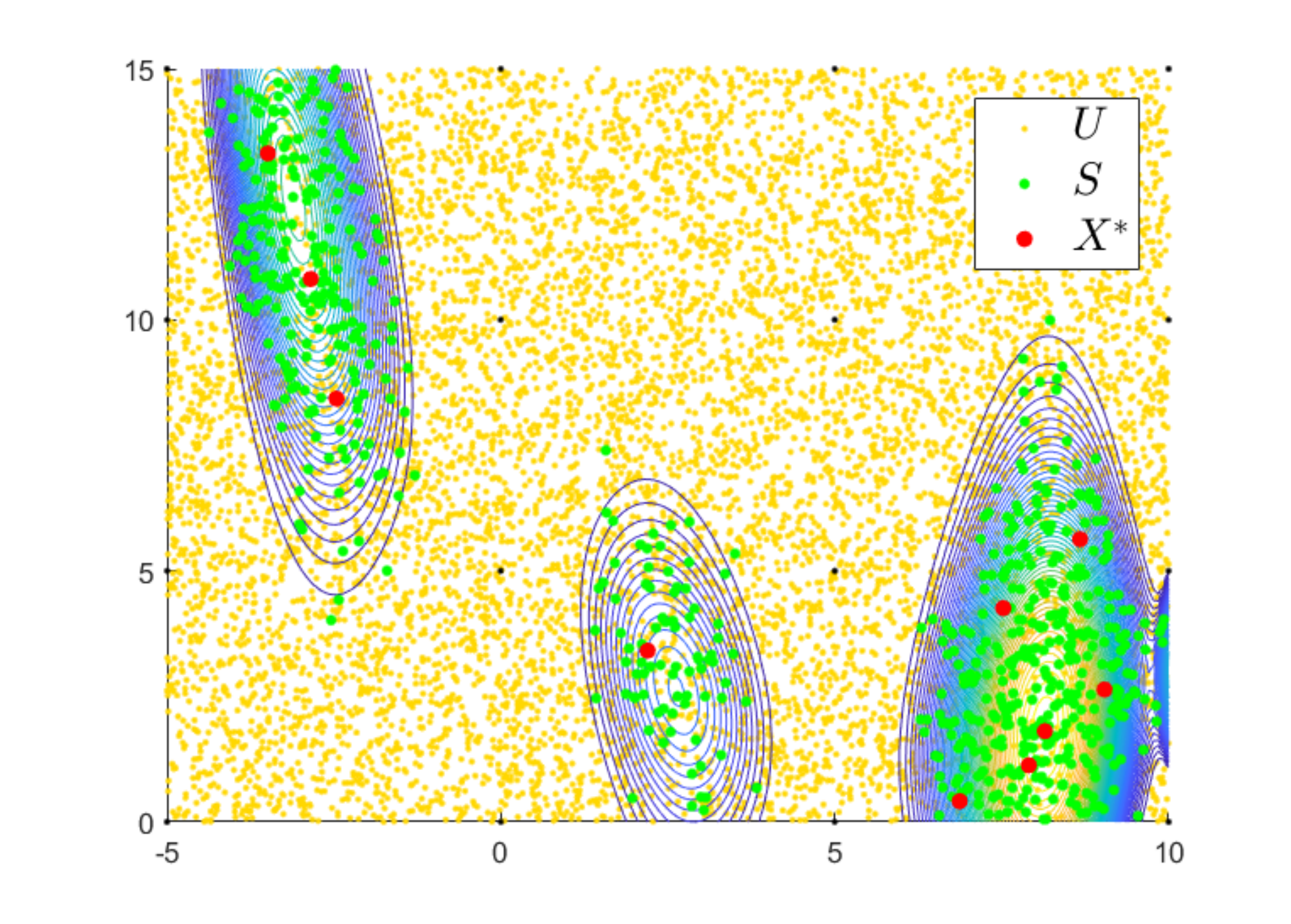}
			\label{fig: exa SCO}
		\end{minipage}%
	}
	\centering
	\subfigure[Uncertainty of SCO and sampling-only]{
		\begin{minipage}[t]{.45\textwidth}
			\centering
			\includegraphics[width=\textwidth]{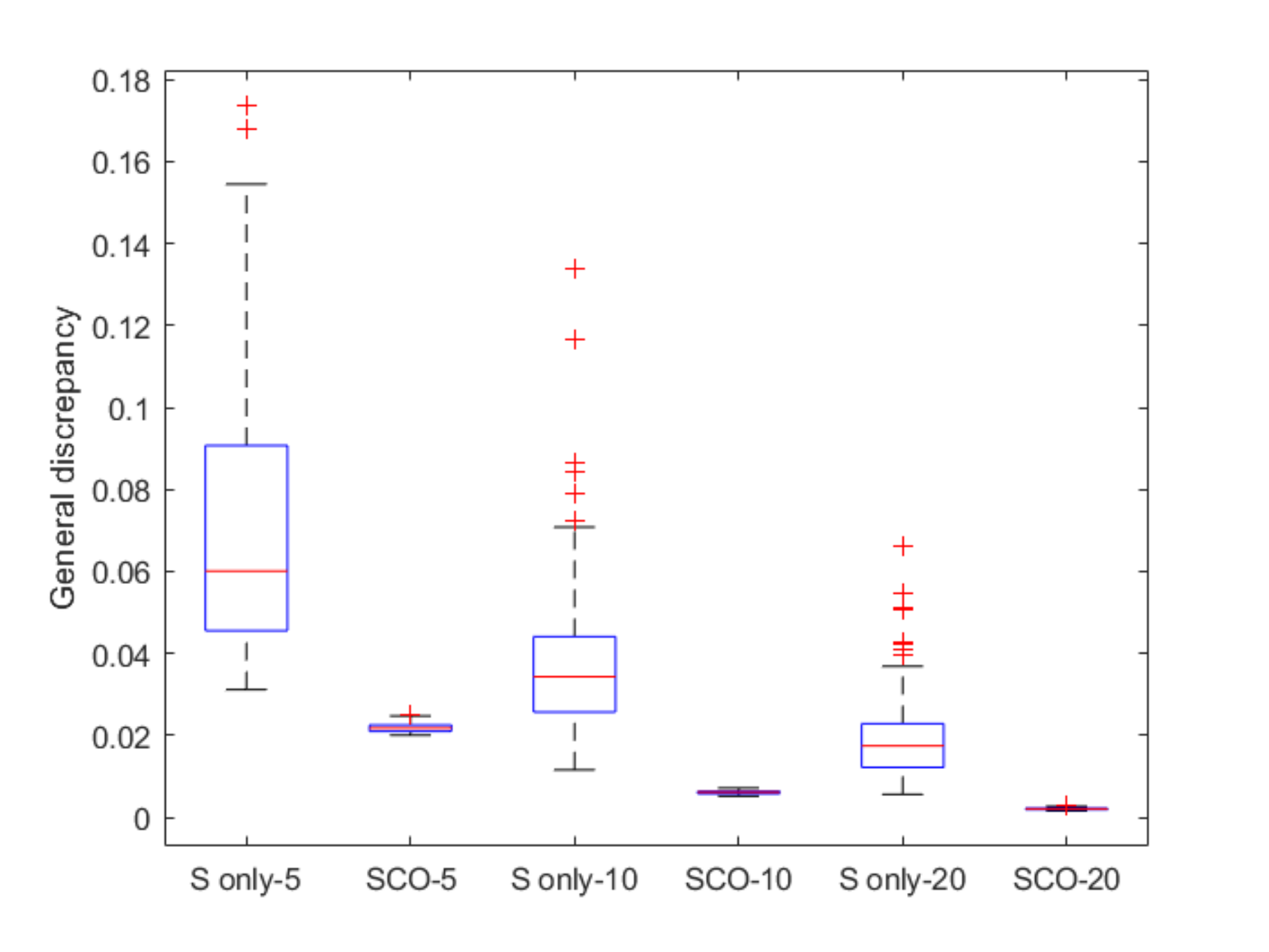}
			\label{fig: compare uncertainty}
		\end{minipage}%
	}
	\caption{SCO example}
\end{figure}

The flowchart of the SCO is shown in Figure \ref{fig: flow SCO}, and we illustrate its process with an example.
We build a Gauss process model according to $Y_{\rm data}=f_{\rm exp}(X_{\rm data})$, where $f_{\rm exp}$ is the Branin function
(\href{http://www.sfu.ca/~ssurjano/optimization.html}{www.sfu.ca/\textasciitilde ssurjano/optimization.html}),
and $X_{\rm data}$ is the set of $4\times4$ mesh points on $\Omega=[-5,10]\times[0,15]$.
We then take EI as the acquisition function \citep{jones1998efficient} and draw the contours of EI, see figure \ref{fig: exa SCO}.
The EI criterion captures the three extreme points of Branin function well.
Then the pre-sample set $U$, sample set $S$, and the final design $X^*$ are plotted in different color.
The size of $U$ is large, and it does not contain the features of $\phi$.
Compared with the optimization from $U$, the optimization from $S$ can greatly reduce the number of possible combinations without losing the features of $\phi$, and $X^*$ is less uncertain than $S$.
This point of view can be seen in Figure \ref{fig: compare uncertainty}.
We repeated 100 tests to compare the uncertainty of the two methods---the sampling-only method and the SCO method.
We chose design sizes of 5,10, and 20 and compare their general discrepancies with the boxplots, which were grouped as S only-5, SCO-5 and so on.
From Figure \ref{fig: compare uncertainty}, we can see the uncertainty of the SCO designs were much less than the uncertainty of the sampling-only designs.
The SCO designs were even better than sampling-only designs with the double sample size, in terms of general discrepancy.

\section{Numerical results}\label{sec num}

In this section, we first study the efficiency of the GA and SA in the SCO.
We then compare some batch Bayesian optimization methods from different perspectives.

\subsection{GA and SA in SCO}

\begin{figure}[htb]
	\centering
	\subfigure{
		\begin{minipage}[t]{\textwidth}
			\centering
			\includegraphics[width=\textwidth]{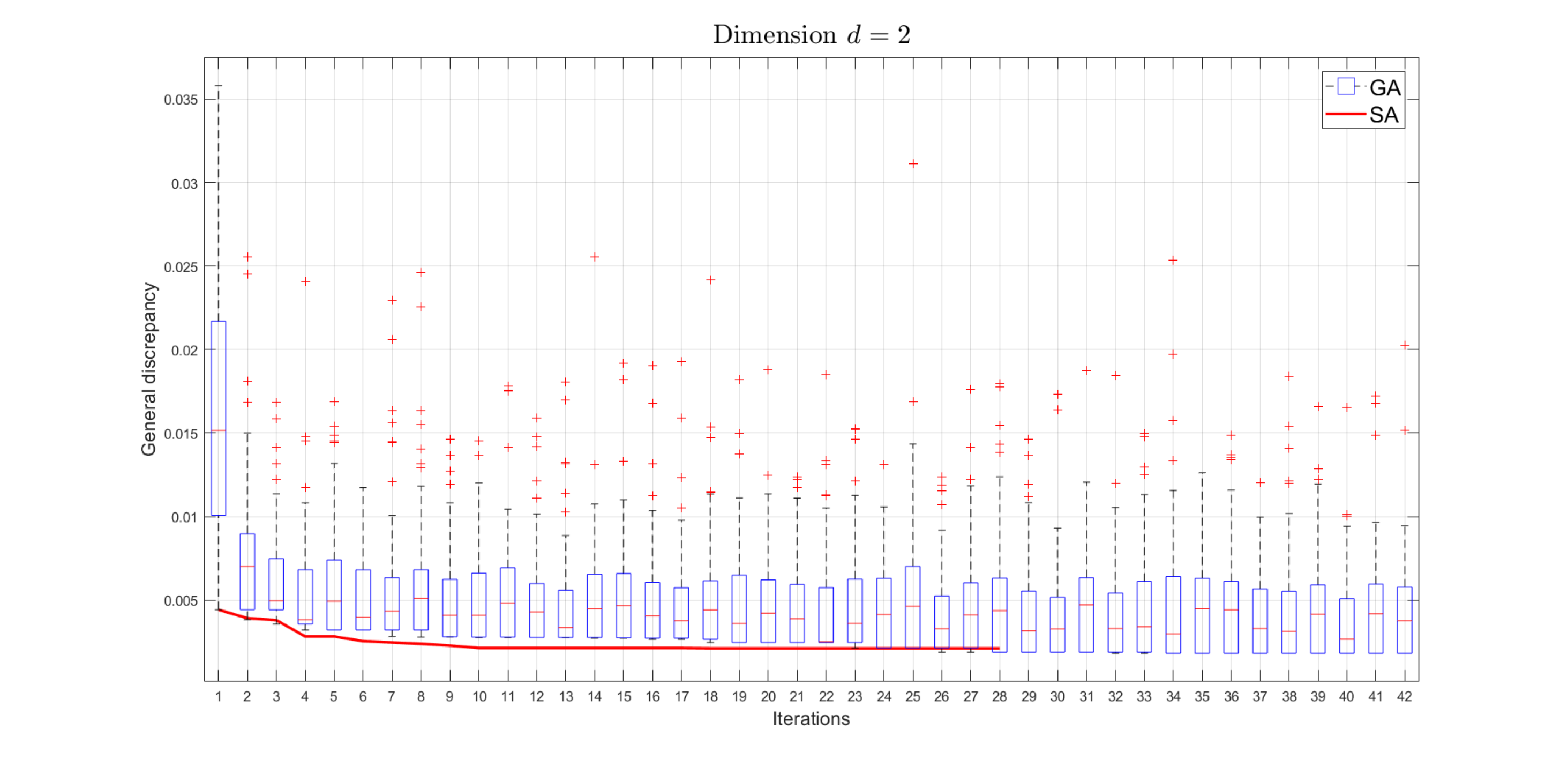}
		\end{minipage}%
	}
	\centering
	\subfigure{
		\begin{minipage}[t]{\textwidth}
			\centering
			\includegraphics[width=\textwidth]{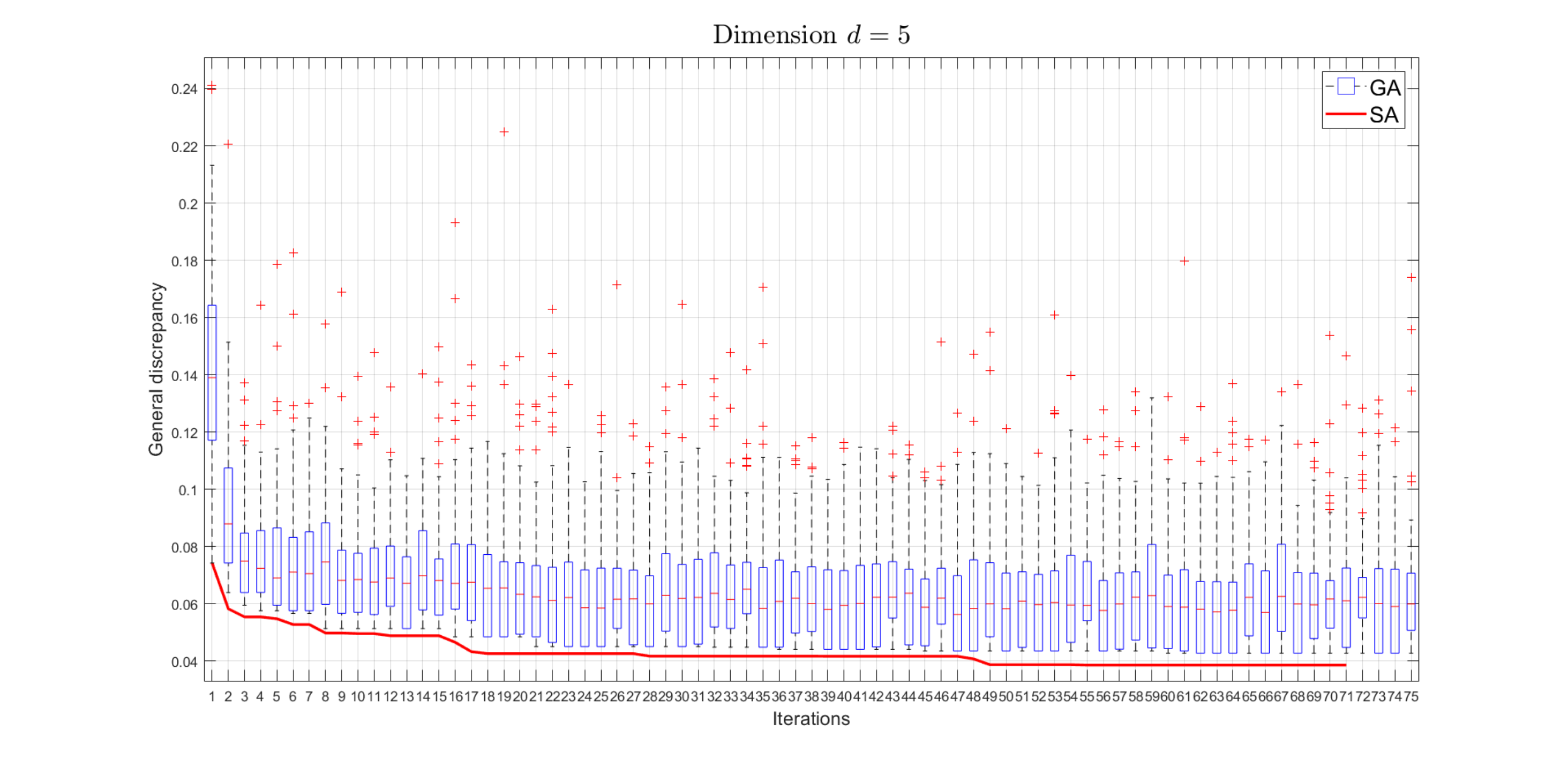}
		\end{minipage}%
	}
\caption{Comparison of the GA and SA in different dimensions}
\label{fig: co GS1}
\end{figure}

In Section \ref{SCO-O}, we introduce the GA and SA for optimization.
In order to compare their efficiencies, we generated the same inputs by Algorithm \ref{SC}, then proceeded to optimization by Algorithms \ref{GA} and \ref{SA}, respectively.
The general discrepancies of the two algorithms are shown in Figure \ref{fig: co GS1}.
The boxplots record the evolution of the population in the GA, and the red line records the switching design in the SA.
From Figure \ref{fig: co GS1}, we can see that the SA converged faster than the GA, and the GA was able to provide slightly better results when the dimension was lower.
In addition, the SA is monotonically decreasing within the iteration in theory; however, the GA is stochastic.
Although the 4-th line in Algorithm \ref{GA} was adopted to ensure the monotonicity of the best individual, the population of the GA fluctuated greatly during evolution.
In consideration of efficiency and stability, we chose the SA in the optimization of the SCO.

\subsection{Efficiency of batch Bayesian optimization}

As shown in the last section, the SCO design has good properties in terms of general discrepancy.
However, can the advantage of general discrepancy be reflected in the batch Bayesian optimization?
In this subsection, we will answer this question with numerical results.
We compared the SCO method with some existing batch methods, including the $q$-EI method with kriging believer (KB) and constant liar minimum (CL-min) as heuristic strategies \citep{ginsbourger2010kriging}, multi-point sampling based on kriging (MPSK) method \citep{cai2017a} and accelerated EGO (aEGO) method \citep{ning2020batch}.

\begin{figure}[htbp]
	\centering
	\includegraphics[scale=0.6]{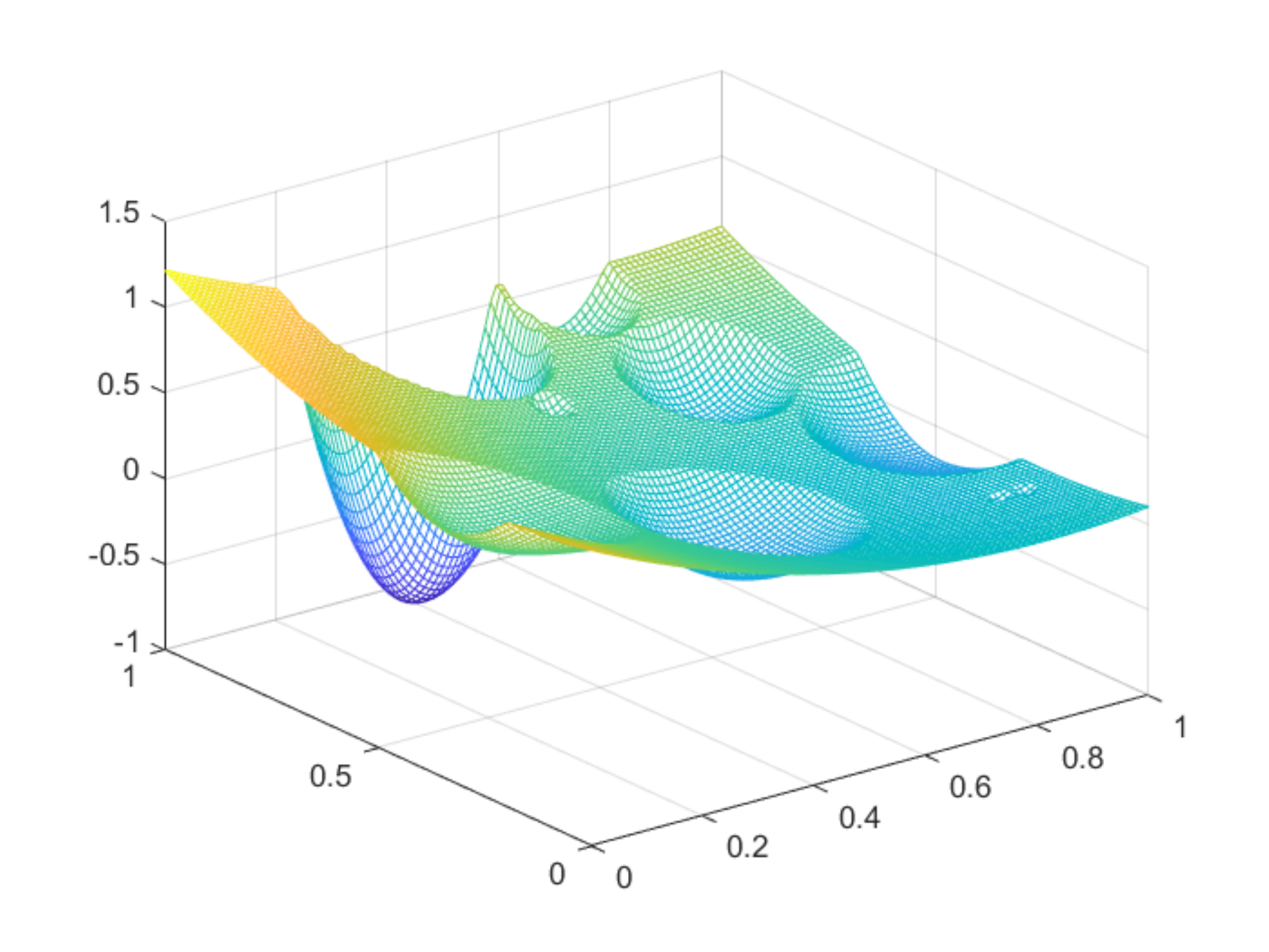}
	\caption{Demonstration of GKLS}
	\label{fig: GKLS}
\end{figure}

We compare SCO method with some existing batch methods, including $q$-EI method with “KB” and “CL-min” as heuristic strategies \citep{ginsbourger2010kriging}, MPSK method \citep{cai2017a} and aEGO \citep{ning2020batch}.
Note that all these methods build a Gauss process model as a surrogate.
However, KB and CL-min are based on a multi-point acquisition function; MPSK, aEGO, and SCO proposed in this article are sample-based methods.
In order to reduce the contingency of the results, we used the Gaviano-Kvasov-Lera-Sergeyev (GKLS) method \citep{1998Test,gaviano2011software} to generate several random functions as the objective.
The parameters were set according to \cite{gaviano2011software}; one of the GKLS functions in two dimensions is demonstrated in Figure \ref{fig: GKLS}.
To compare the efficiency of the optimization, we define the average relative accuracy (ARA) as
\begin{equation}
{\rm ARA}=\frac{1}{N_{\rm f}}\sum_{i=1}^{N_{\rm f}}\frac{y_{\rm min}^{(i)}-f_{\rm min}^{(i)}}{|f_{\rm min}^{(i)}|},
\end{equation}
where $N_{\rm f}$ is the number of random functions, $f_{\rm min}^{(i)}$ is the optimal value of the $i$-th function, and $y_{\rm min}^{(i)}$ is the minimum response found by the algorithm.

\begin{figure}[htbp]
	\centering
	\subfigure{
		\begin{minipage}[t]{.45\textwidth}
			\centering
			\includegraphics[width=\textwidth]{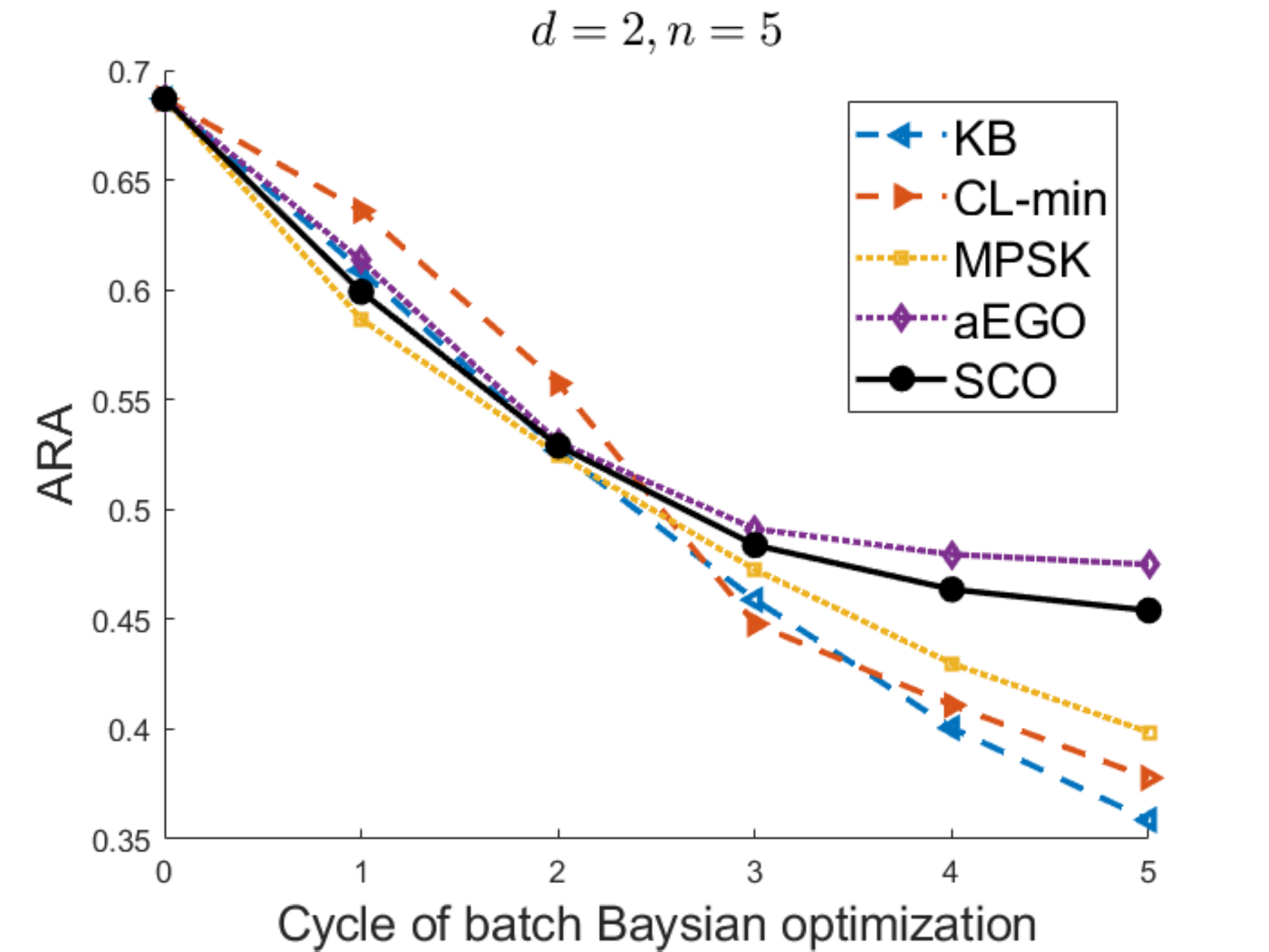}
		\end{minipage}%
	}
	\centering
	\subfigure{
		\begin{minipage}[t]{.45\textwidth}
			\centering
			\includegraphics[width=\textwidth]{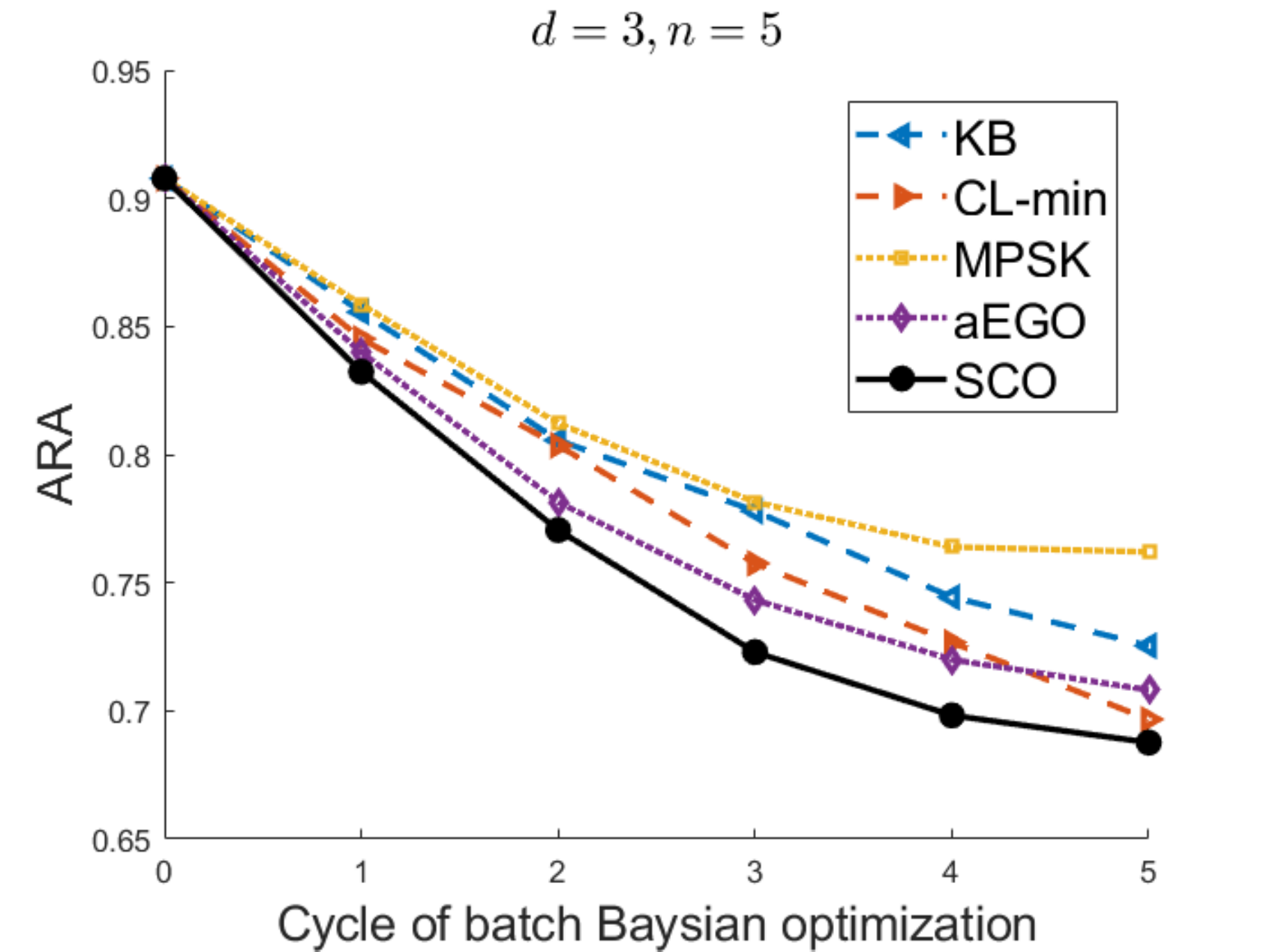}
		\end{minipage}%
	}
\centering
\subfigure{
\begin{minipage}[t]{.45\textwidth}
	\centering
	\includegraphics[width=\textwidth]{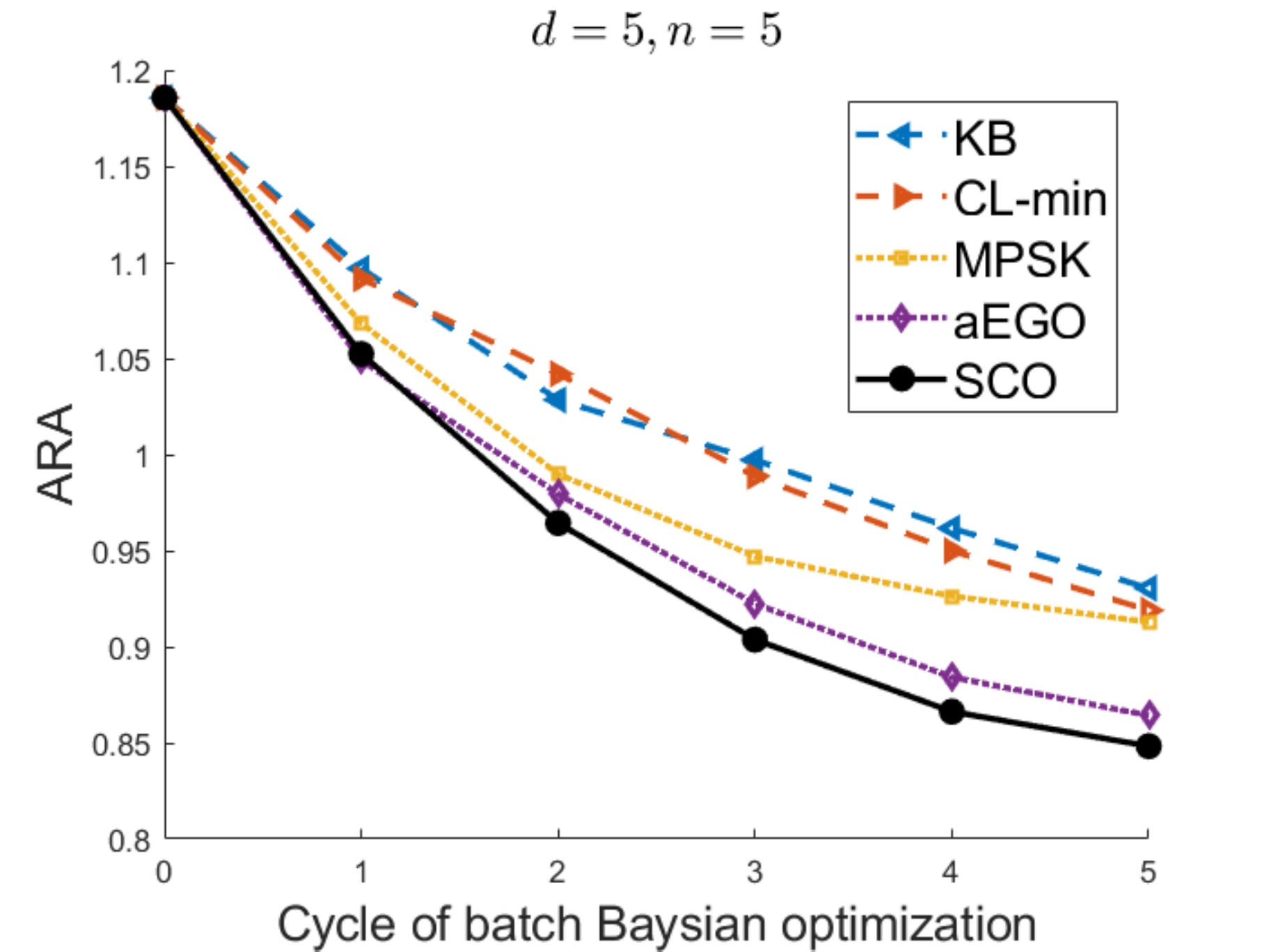}
\end{minipage}%
}
\centering
\subfigure{
\begin{minipage}[t]{.45\textwidth}
	\centering
	\includegraphics[width=\textwidth]{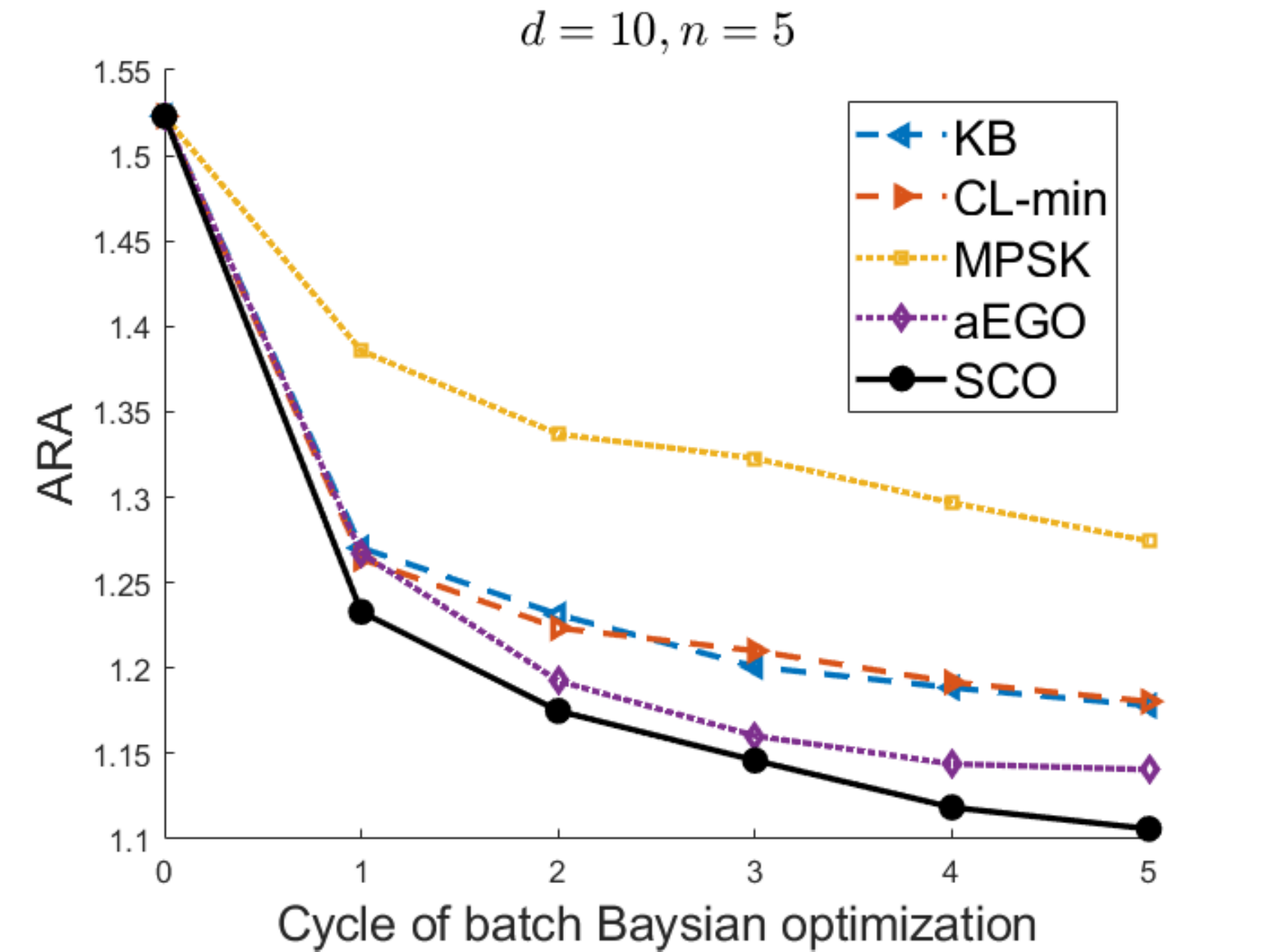}
\end{minipage}%
}
	\caption{Comparison of batch Bayesian optimizations with different dimensions}
	\label{fig: co dimen}
\end{figure}

First, we analyze the efficiency of the batch methods in different dimensions.
We generated 100 functions for dimension $d=2,3,5,\ {\rm and}\ 10$.
For each dimension, we started with same initial data sites $X_{\rm data}$ to optimize them.
The number of $X_{\rm data}$ was set as $5\times d$, and $X_{\rm data}$ were constructed by uniform designs \citep{fang2018theory}.
We conducted five cycles of batch Bayesian optimization, with each batch size set as $n=5$.
Figure \ref{fig: co dimen} displays the ARAs of different batch methods in different dimensions.
As to the methods based on a multi-point acquisition function, KB and CL-min performed well only in two dimensions.
SCO and aEGO showed an advantage in higher dimensions, and SCO performed the best of all the sample-based methods.

\begin{figure}[htbp]
	\centering
	\subfigure{
		\begin{minipage}[t]{.45\textwidth}
			\centering
			\includegraphics[width=\textwidth]{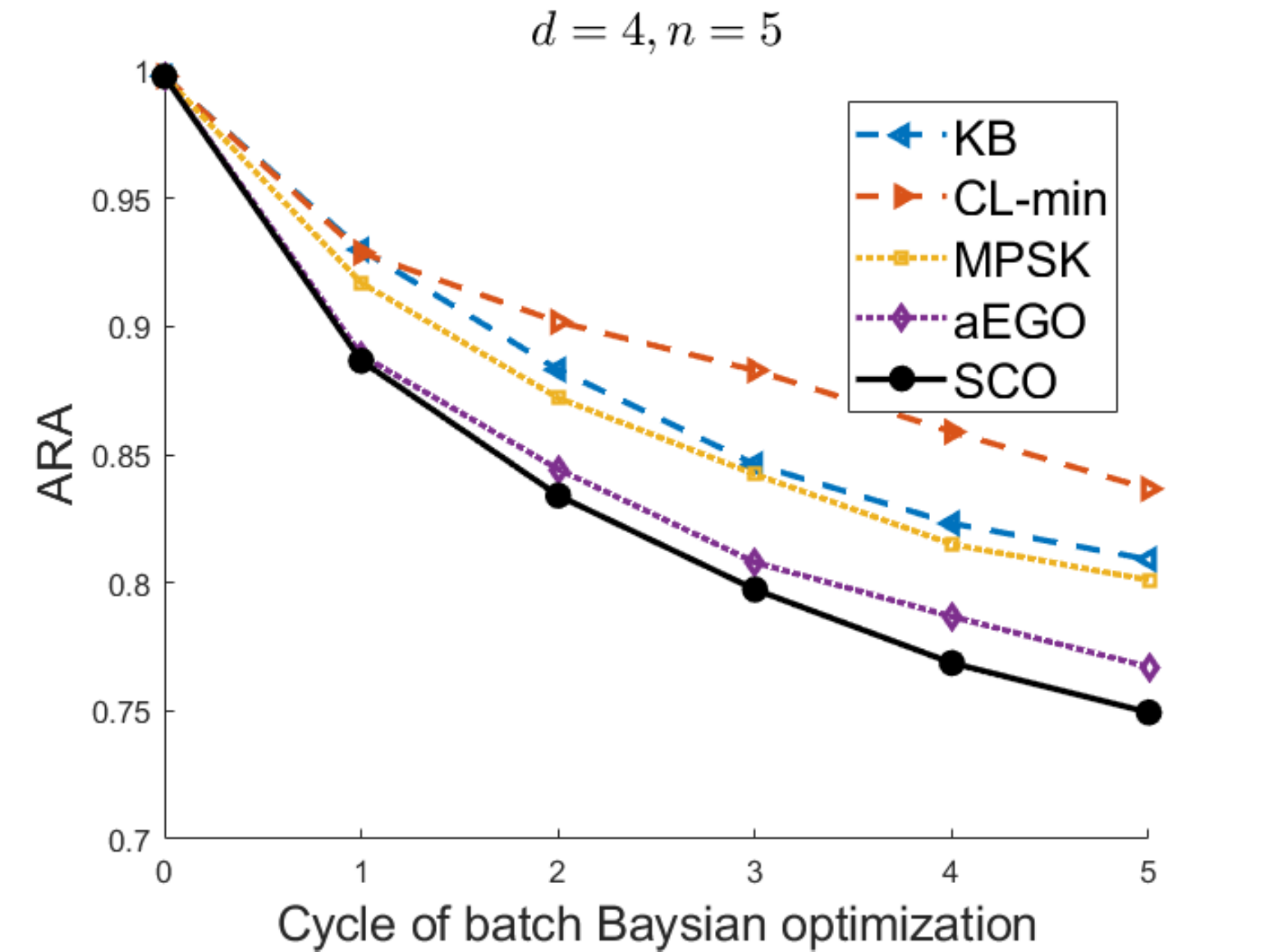}
		\end{minipage}%
	}
	\centering
	\subfigure{
		\begin{minipage}[t]{.45\textwidth}
			\centering
			\includegraphics[width=\textwidth]{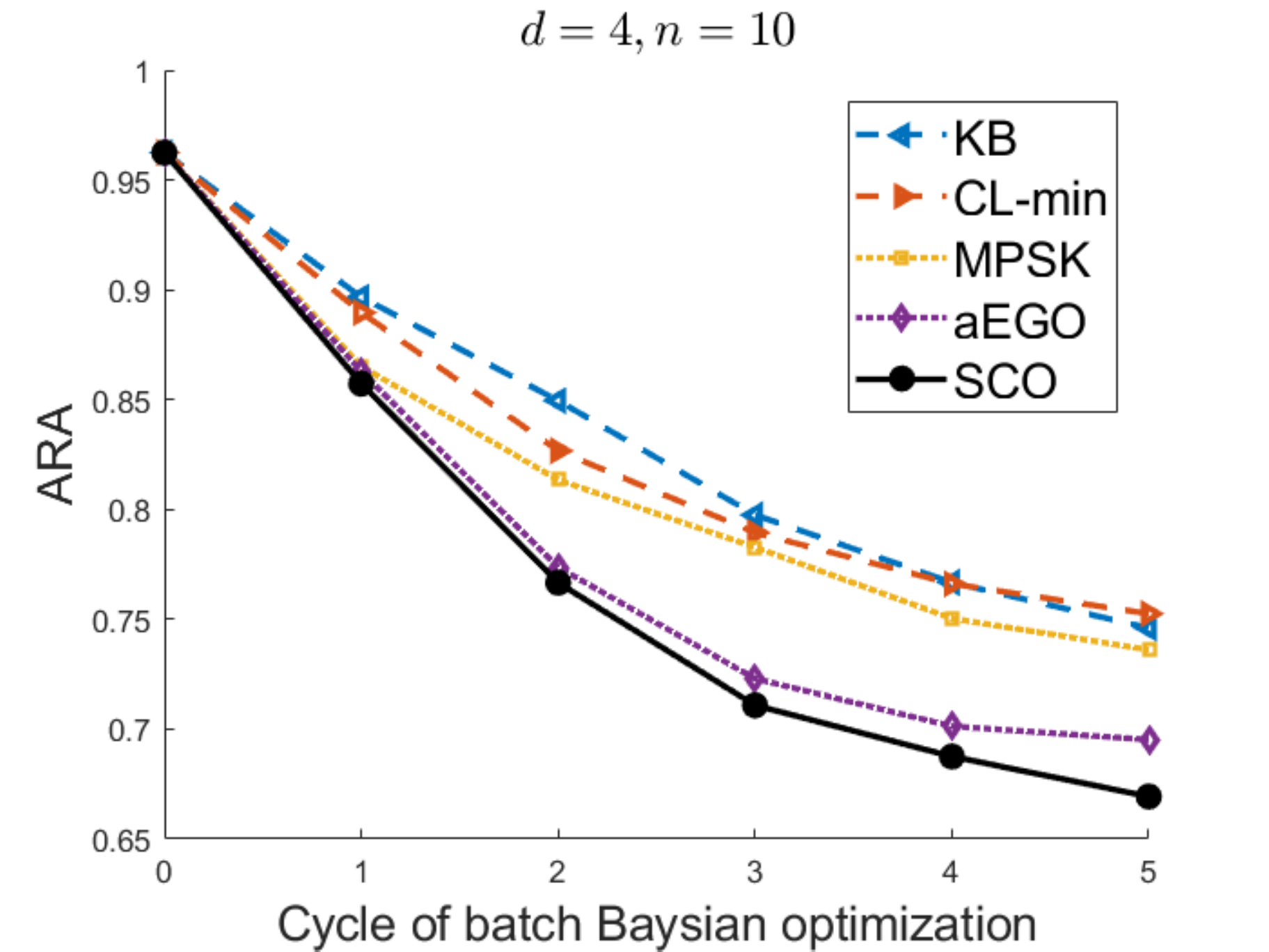}
		\end{minipage}%
	}
	\centering
	\subfigure{
		\begin{minipage}[t]{.45\textwidth}
			\centering
			\includegraphics[width=\textwidth]{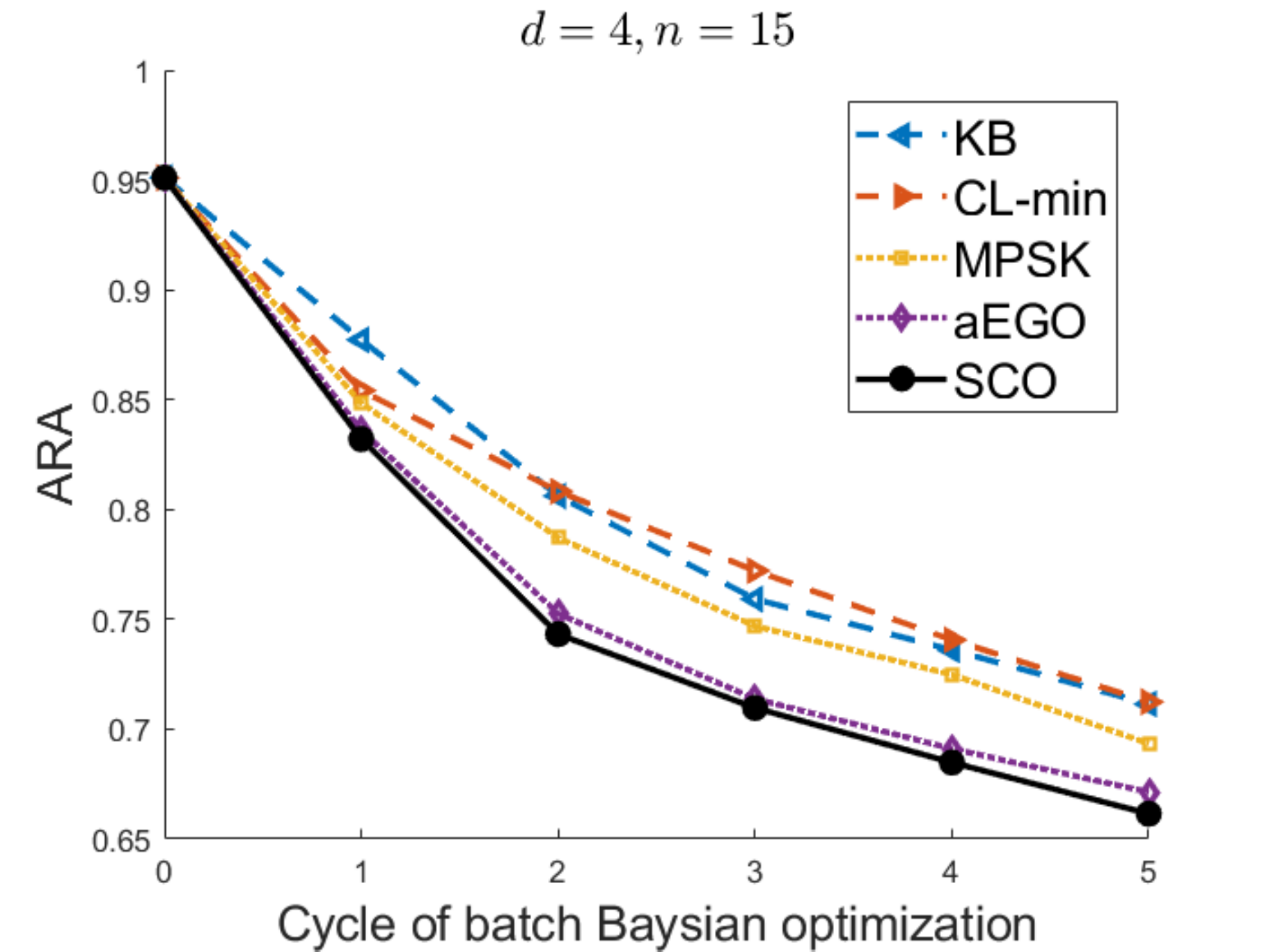}
		\end{minipage}%
	}
\centering
\subfigure{
\begin{minipage}[t]{.45\textwidth}
	\centering
	\includegraphics[width=\textwidth]{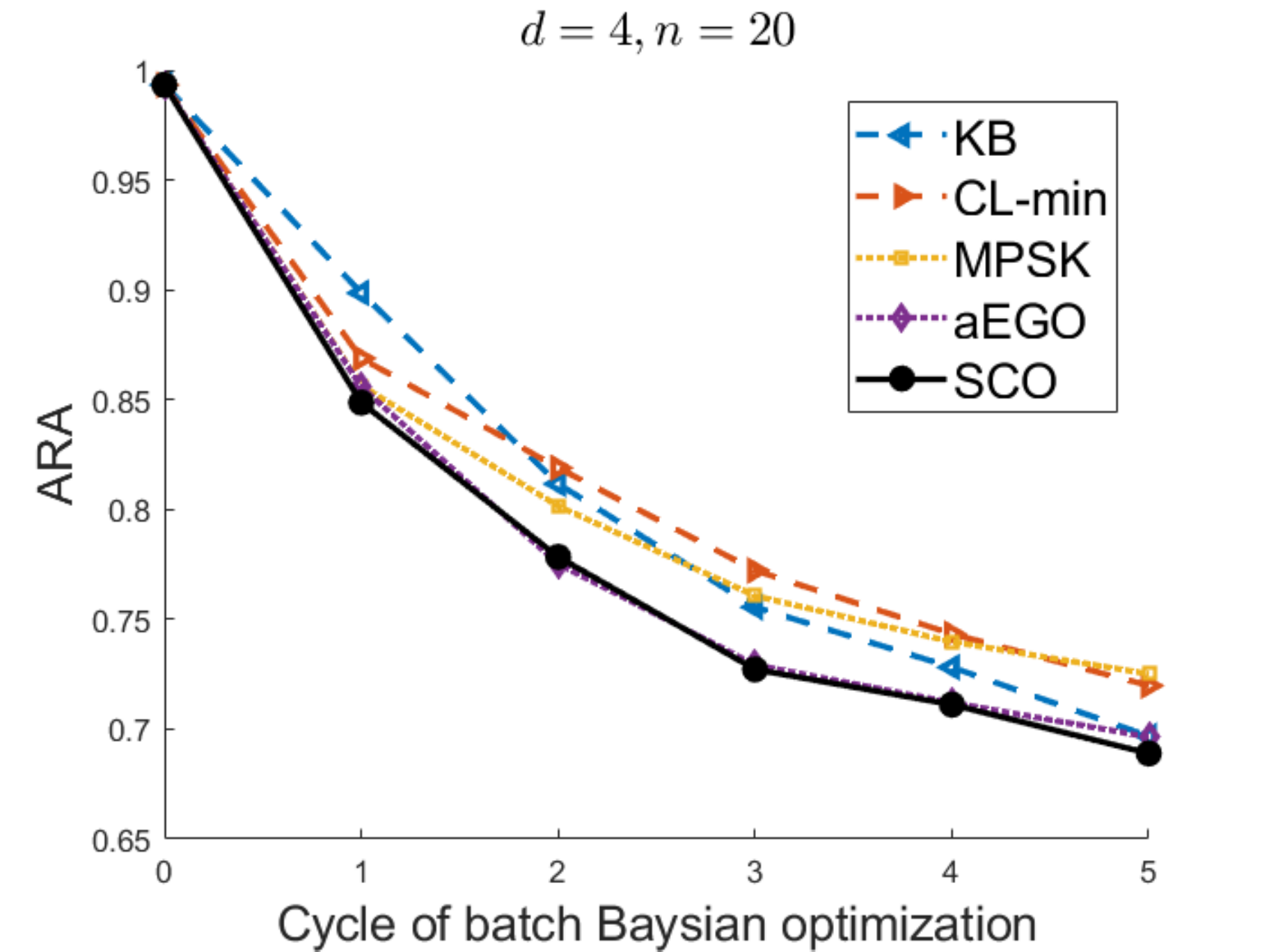}
\end{minipage}%
}
	\centering
	\caption{Comparison of  batch Bayesian optimizations with different batch size}
	\label{fig: co size}
\end{figure}

Second, we analyzed the influence of batch size.
For $d=4$, we conducted similar contrast experiments, but set $n=5,10,15,\ {\rm and}\ 20$.
The ARAs are displayed in Figure \ref{fig: co size}.
As the batch size increased, the difference between SCO and aEGO decreases, the ARAs of SCO and aEGO almost overlap when $n=20$.
This illustrates that, the effect of calculation and optimization is not obvious when the batch size is large.
It stands to reason that if we have enough experimental resources, experimental design would not matter.
However, SCO is still most efficient among all sampling-based methods.

\subsection{Application in radar interference experiment}\label{sec exa}

\begin{figure}[htbp]
	\centering
	\includegraphics[scale=0.6]{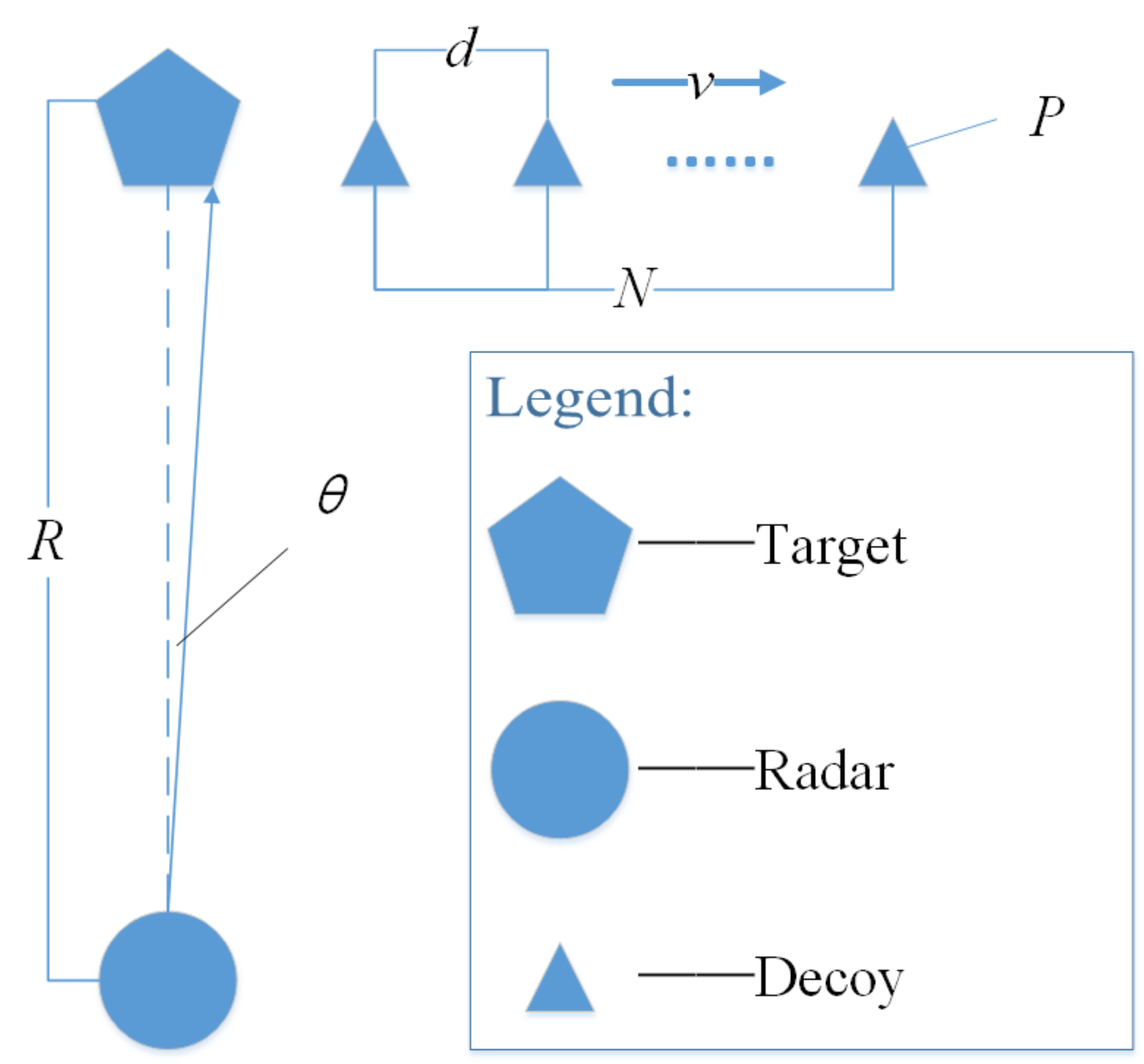}
	\caption{Schematic diagram of single experimental scene}
	\label{fig: single experimental scene}
\end{figure}

The interference ability test of a certain type of radar is carried out.
The interference method of releasing array decoy is adopted, as shown in Figure \ref{fig: single experimental scene}.
The radar is fixed and the target is set at azimuth 0 and range $R$.
After the radar is turned on, the target releases $N$ decoys spaced $d$, the decoys sail at speed $v$ and release interference signals at power $P$.
The process of target recognition is recorded, and the interference is measured by azimuth angle $\theta$ of radar.
The purpose of the experiment is to find the optimal strategy $(R,d,v,P)$ to maximize $\theta$ when $N$ is fixed.

\begin{figure}[htbp]
	\centering
	\includegraphics[scale=0.6]{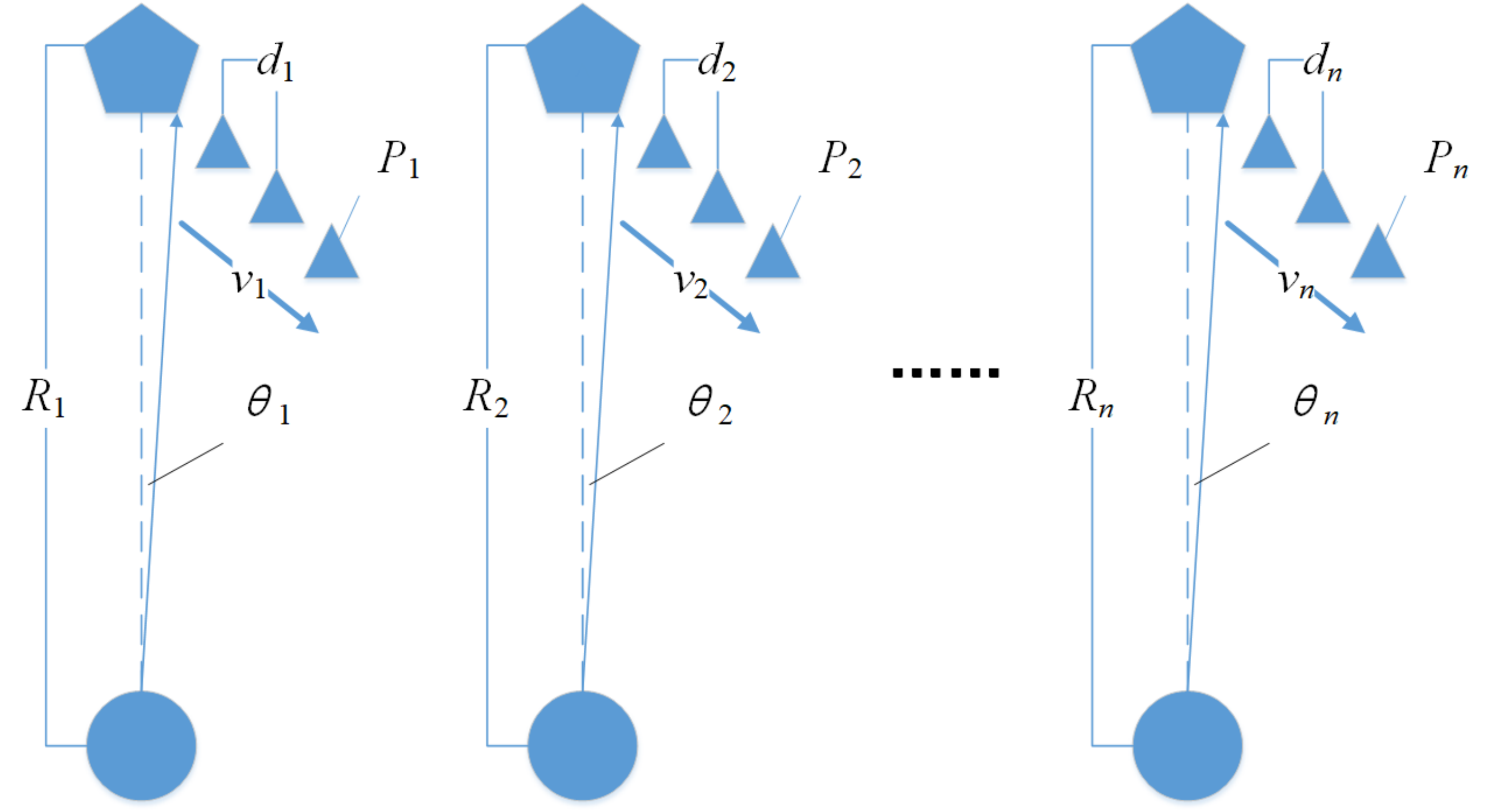}
	\caption{Schematic diagram of batch experimental scene}
	\label{fig: batch experimental scene}
\end{figure}

\begin{figure}[htbp]
	\centering
	\includegraphics[width=\textwidth]{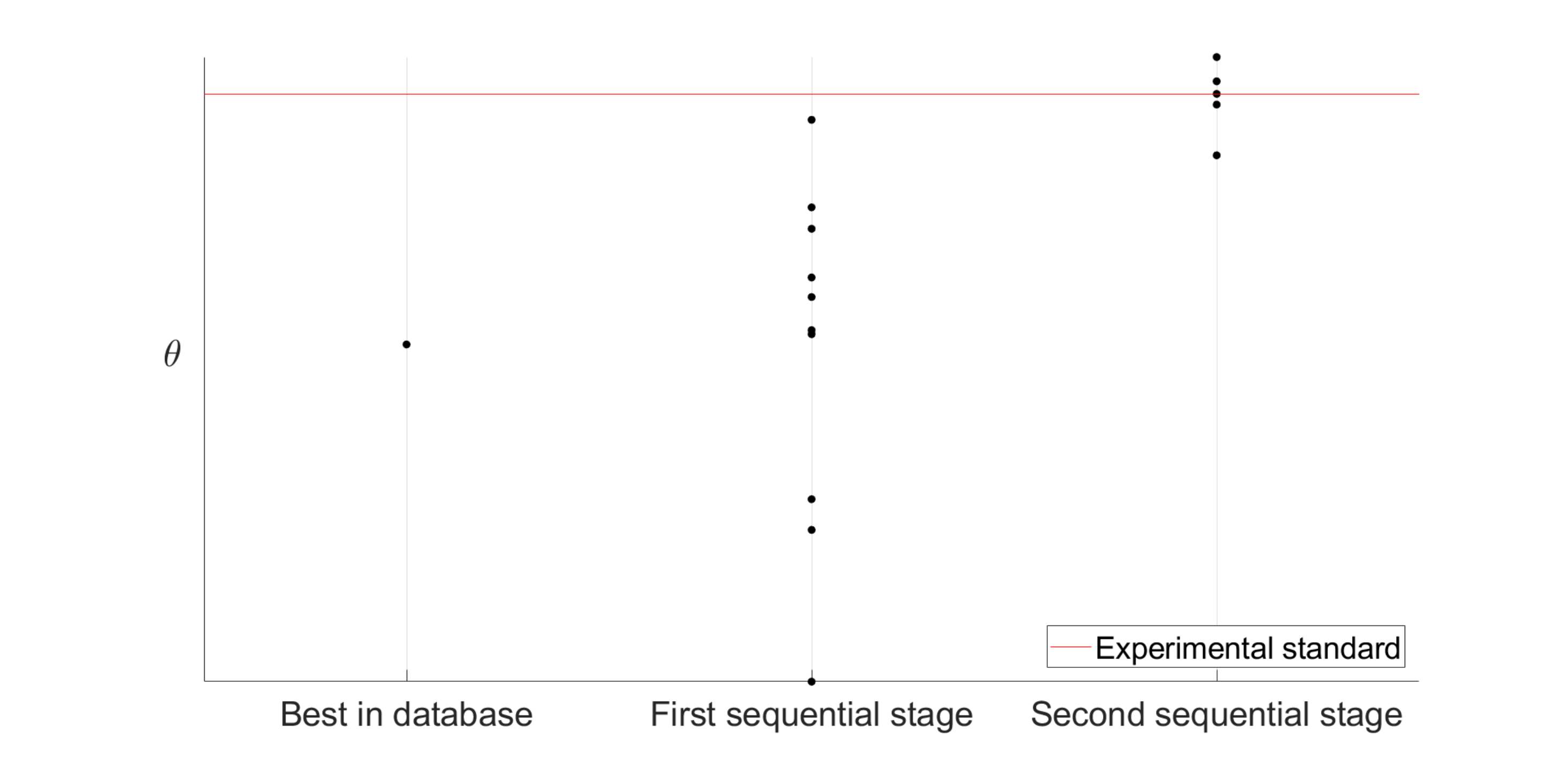}
	\caption{Experimental results}
	\label{fig: Experimental results}
\end{figure}

Preliminary modeling is conducted based on previous experimental data and sequential experiments are carried out in the next stage.
$n$ experiments can be carried out simultaneously in the same stage, as shown in Figure \ref{fig: batch experimental scene}.
In order to shorten the experimental period, the batch sequential experiments were carried out, and the SCO method was used to design the batch sequential experiments.
Considering the resources of the experimental base, $n_1$ experiments were designed in the first stage.
Analysis of the results showed that the experimental standard was still not met and further experiments were needed.
Due to the change of the resources in the experimental base, $n_2$ sequential experiments were designed in the second stage.
The SCO has the flexibility to deal with different batch size.
The experimental results of the two stages are shown in Figure \ref{fig: Experimental results}, where the interference $\theta$ is improved.
Analyses confirmed that the results met the experimental standard.

In this case, we made full use of experimental resources and met the experimental standard in two stages of batch experiments, which greatly shortened the experimental period.
The SCO method is used for batch sequential experimental design, which is flexible and less uncertain.
This means the SCO is more suitable for complex and high-cost experiments.

\section{Conclusion and further discussion}

This article introduced a sequential design method for batch Bayesian optimization.
The main processes of this method include sampling-calculation-optimization (SCO).
SCO is a sampling-based method that does not construct a new acquisition function but samples from the existing acquisition function.
To reduce the uncertainty, the samples are optimized in the sense of general discrepancy.
We have proposed several strategies to reduce the amount of computation and make the optimization possible.
Numerical results show that the uncertainty of the SCO was much less than the sampling-only method.
In addition, the batch Bayesian optimization with SCO was more efficient than other batch methods.
Although we have introduced the SCO based on the Gaussian process model and the EI criterion, the method is also well-suited to other models and acquisition functions.
Finally, the case of radar interference experiment shows the application value and scenario of SCO method.

In this article, the batch size is fixed and constrained by the experimental resource of every experimental period.
In other kind of experiments, the batch size may be flexible and it is a part of the experimental design.
In such case, how to design the batch size to experiment more efficiently remains to be further studied.

\section*{Acknowledgments}

This work is supported by the National Natural Science Foundation of China (No. 11771450, 12101608).

%

\bibliographystyle{apalike}
\bibliography{SCO_arXiv}

\end{document}